\documentclass[10pt,titlepage,oneside]{article}
\usepackage[dvips]{graphics}
\usepackage{epsfig,rotating}
\usepackage{amssymb,amsmath}
\usepackage{latexsym}

\newcommand{\be}{\begin{equation}}
\newcommand{\ee}{\end{equation}}
\newcommand{\beaa}{\begin{eqnarray*}}
\newcommand{\eeaa}{\end{eqnarray*}}
\newcommand{\bea}{\begin{eqnarray}}
\newcommand{\eea}{\end{eqnarray}}
\newcommand{\lbl}{\label}

\newcommand{\ml}{\mathcal}

\newcommand{\bd}{\bold}

\newtheorem{theorem}{ \noindent T{\footnotesize HEOREM}}

\newtheorem{lemma}{ \noindent L{\footnotesize EMMA}}[section]

\newtheorem{remark}{ \noindent R{\footnotesize EMARK}}[section]

\newtheorem{theoremm}{ \noindent T{\footnotesize HEOREM} A.\!\!}

\begin{document}

\begin{centering}
{\large \bf Limit Theorems for Beta-Jacobi Ensembles}
\\[2em]

{\sc Tiefeng Jiang} \noindent\footnote{Supported in part by
NSF\#DMS-0449365, School of Statistics, University of Minnesota, 224
Church Street, MN55455, tjiang@stat.umn.edu.
\newline AMS (2000) subject classifications: 60B10, 60B20, 60F05, 60F15, 62H10.}


\vskip.2in

\end{centering}

\noindent{\bf Abstract}\ For a $\beta$-Jacobi ensemble determined by parameters $a_1, a_2$ and $n$, under the restriction that the three parameters go to infinity with $n$ and $a_1$ being of small orders of $a_2$, we obtain both the bulk and the edge scaling limits. In particular, we derive the asymptotic distributions for the largest  and the smallest eigenvalues, the Central Limit Theorems of the eigenvalues, and the limiting distributions of the empirical distributions of the eigenvalues.
\vspace{.1cm}\\

\section{Introduction}\lbl{intro}
\setcounter{equation}{0}
Let $\beta>0$ be a constant and $n\geq 2$ be an integer. A beta-Jacobi ensemble, also called in the literature as the
beta-MANOVA ensemble, is a set of random variables $(\lambda_1,\cdots, \lambda_n)\in [0,1]^n$ with probability density function
\begin{eqnarray}\lbl{beJacobi}
f_{\beta, a_1, a_2}(\lambda)=c_{J}^{\beta, a_1, a_2}\prod_{1\leq
i<j\leq
n}|\lambda_i-\lambda_j|^{\beta}\cdot\prod_{i=1}^n\lambda_i^{a_1-p}(1-\lambda_i)^{a_2-p},
\end{eqnarray}
 where $a_1,\, a_2>\frac{\beta}{2}(n-1)$ are parameters, $p=1+ \frac{\beta}{2}(n-1)$,  and
\begin{eqnarray}\lbl{conJacobi}
c_{J}^{\beta, a_1,
a_2}=\prod_{j=1}^{n}\frac{\Gamma(1+\frac{\beta}{2})\Gamma(a_1+a_2-\frac{\beta}{2}(n-j))}
{\Gamma(1+\frac{\beta}{2}j)\Gamma(a_1-\frac{\beta}{2}(n-j))\Gamma(a_2-\frac{\beta}{2}(n-j))}.
\end{eqnarray}
The ensemble is associated with the multivariate analysis of
variance (MANOVA). For $\beta=1,2$ and $4,$ the function
$f_{\beta}(\lambda)$ in (\ref{beJacobi}) is the density function of
the eigenvalues of $\bd{Y}^*\bd{Y}(\bd{Y}^*\bd{Y}+
\bd{Z}^*\bd{Z})^{-1}$ with $a_1=\frac{\beta}{2}m_1$ and
$a_2=\frac{\beta}{2}m_2$, where $\bd{Y}=\bd{Y}_{m_1\times n}$ and $\bd{Z}=\bd{Z}_{m_2\times n}$ are independent matrices with $m_1, m_2\geq n$, and the entries of both
 matrices are independent random variables with the standard real, complex or quaternion
 Gaussian distributions. See \cite{Cons} and \cite{Muirhead1982} for
 $\beta=1,2.$ Other references about the connections between the Jacobi ensembles and statistics are \cite{Anderson, Collins, Collins03, Cons, DEL, Eaton, Conductance, Jiang07, Jiang06, Muirhead1982}.

In statistical mechanics, the model of the log gases can be characterized by the beta-Jacobi ensembles.
A log gas is a system of charged particles on the real line which are subject to
a logarithmic interaction potential and Brownian-like changes. If
the particles are contained in the interval [0, 1] and are also subject to the
external potential $\sum_{i=1}^n(\frac{r+1}{2}-\frac{1}{\beta})\log \lambda_i + (\frac{s+1}{2}-\frac{1}{\beta})\log (1-\lambda_i)$, where $r=\frac{2}{\beta}a_1 -n$ and $s=\frac{2}{\beta}a_2 -n$, and $\beta$ is the inverse of the temperature,
then it is known that the stationary distribution of the system of charges in the long term is the Jacobi
ensemble as in (\ref{beJacobi}), see, e.g., \cite{Baker, Dyson, Forresterbook, Wang}.

The beta-Jacobi ensembles also have connections to other subjects in  mathematics and
physics, for instance, lattice gas theory \cite{Forresterbook, Conductance},  Selberg integrals \cite{Forrester, Macdonald1998, Okounkov} and Jack functions \cite{Aomoto, Kadell, stanley}.

Now we briefly recall some research on the beta-Jacobi ensembles. Lippert \cite{Lippert} gives a model to generate the beta-Jacobi ensembles (see also \cite{Nenciu} for a similar method used in the construction of the beta-circular ensembles). In studying the largest principal angles between random subspaces,
Absil, Edelman and Koev \cite{Absil} obtain a formula related to the Jacobi ensembles. Edelman and Sutton \cite{sutton} study CS decomposition and singular values about these models. Dumitriu and Koev \cite{Dumitriu2008} derive the exact distributions of the largest eigenvalues for the ensembles. Jiang \cite{Jiang07} derives the bulk and the edge scaling limits for the beta-Jacobi ensembles for  $\beta=1, 2$ when $p$ and $a_1$ in (\ref{beJacobi}) are of small orders of $a_2$. Johnstone \cite{Johnstone} obtains the asymptotic distribution of the largest eigenvalues  for $\beta=1, 2$ when $a_1, a_2$ and $p$ in (\ref{beJacobi}) are proportional to each other. Recently, Demni \cite{Demni} investigates the beta-Jacobi processes.

In this paper, for the beta-Jacobi ensembles, we study the asymptotic distributions of the largest and smallest eigenvalues, the limiting empirical distributions of the eigenvalues, the law of large numbers and the central limit theorems for the eigenvalues. Before stating the main results, we need some notation.

Let $\beta>0$ be a fixed constant, $n\geq 2$ be an integer, $a_1$ and $a_2$ be positive variables.  The following condition will be used later.
\begin{eqnarray}\lbl{shanxi}
n\to \infty,\ a_1\to \infty\ \mbox{and}\ a_2\to \infty\ \mbox{such that}\ a_1=o(\sqrt{a_2}),\ n=o(\sqrt{a_2})\ \mbox{and}\ \frac{n\beta}{2a_1}\to \gamma\in (0, 1].
\end{eqnarray}
For two Borel probability measures $\mu$ and $\nu$ on $\mathbb{R}^k$, recall the metric
\begin{eqnarray}\lbl{du}
d(\mu, \nu)=\sup_{\|f\|_{BL}\leq 1}\left|\int_{\mathbb{R}^k} f(x)\,d\mu - \int_{\mathbb{R}^k} f(x)\,d\nu\right|,
\end{eqnarray}
where $f(x)$ is a bounded Lipschitz function defined on $\mathbb{R}^k$ with
\begin{eqnarray*}
\|f\|_{BL}=\sup_{x \ne y}\frac{|f(x)- f(y)|}{\|x-y\|} + \sup_{x\in \mathbb{R}^k}|f(x)|.
\end{eqnarray*}
Then, for a sequence of probability measures $\{\mu_n;\, n=0,1,2,\cdots\}$ defined on $(\mathbb{R}^k, \mathcal{B}(\mathbb{R}^k))$, we know $\mu_n$ converges weakly to $\mu_0$ if and only if $d(\mu_n, \mu_0) \to 0$ as $n\to\infty$, see, e.g., \cite{Dudley}. Similarly,  we say that a sequence of random variables $\{Z_n;\, n\geq 1\}$ taking values in $\mathbb{R}^k$ converges  weakly (or in distribution) to a Borel probability measure $\mu$ on $\mathbb{R}^k$  if $Ef(Z_n)\to \int_{\mathbb{R}^k}f(x)\,\mu(dx)$ for any bounded and continuous function $f(x)$ defined on $\mathbb{R}^k.$ This is also  equivalent to that $d(\mathcal{L}(Z_n), \mu) \to 0$ as $n\to \infty,$ where $\mathcal{L}(Z_n)$ is the probability distribution of $Z_n$, see also \cite{Dudley}.

For $\gamma \in (0, 1],$ let $\gamma_{min}=(\sqrt{\gamma}-1)^2$ and
$\gamma_{max}=(\sqrt{\gamma}+1)^2$. The Marchenko-Pastur law is the probability distribution with density function
\begin{eqnarray}\lbl{mooncake}
f_{\gamma}(x)=\begin{cases}
\frac{1}{2\pi \gamma x}\sqrt{(x-\gamma_{min})(\gamma_{max}-x)}\ , & if\ \text{$x\in [\gamma_{min}, \gamma_{max}];$}\\
0, & \text{otherwise.}
\end{cases}
\end{eqnarray}
The following is about the limiting distribution of the empirical eigenvalues of the beta-Jacobi ensembles.

\begin{theorem}\lbl{snowf} Let  $\lambda_1,\cdots,
\lambda_n$ be random variables with density function $f_{\beta, a_1, a_2}(\lambda)$ as in
 (\ref{beJacobi}).  Set
  \begin{eqnarray*}
\mu_n=\frac{1}{n}\sum_{i=1}^n I_{\frac{a_2}{n}\lambda_i}
  \end{eqnarray*}
for $n\geq 2.$ Assuming (\ref{shanxi}), then $d(\mu_n, \mu_0)$ converges to zero in probability, where $\mu_0$ has density $c\cdot f_{\gamma}(c x)$ with $c=2\gamma/\beta$ and $f_{\gamma}(x)$ is as in (\ref{mooncake}).
\end{theorem}
The next result gives the weak laws of large numbers of the largest and smallest eigenvalues for the beta-Jacobi ensembles.

\begin{theorem}\lbl{branch} Let  $\lambda_1,\cdots,
\lambda_n$ be random variables with density function $f_{\beta, a_1, a_2}(\lambda)$ as in
 (\ref{beJacobi}). Set  $\lambda_{max}(n)=\max\{\lambda_1,\cdots, \lambda_n\}$, and $\lambda_{min}(n)=\min\{\lambda_1,\cdots, \lambda_n\}$. Assuming (\ref{shanxi}), we have that
\begin{eqnarray*}
\frac{a_2}{n}\cdot \lambda_{max}(n) \to \beta\cdot \frac{(1+ \sqrt{\gamma}\,)^2}{2\gamma}\ \ \ \mbox{and}\ \ \  \frac{a_2}{n}\cdot \lambda_{min}(n) \to \beta \cdot \frac{(1-\sqrt{\gamma}\,)^2}{2\gamma}
\end{eqnarray*}
in probability.
\end{theorem}
Here is the central limit theorem for the eigenvalues for the model in (\ref{beJacobi}).

\begin{theorem}\lbl{hay} Let  $\lambda_1,\cdots,
\lambda_n$ be random variables with density function $f_{\beta, a_1, a_2}(\lambda)$ as in
 (\ref{beJacobi}).  Given integer $k\geq 1,$ define
\begin{eqnarray*}
X_i=\sum_{j=1}^n\left(\frac{c\,a_2}{n}\lambda_j\right)^i-n
\sum_{r=0}^{i-1}\frac{1}{r+1}\binom{i}{r}\binom{i-1}{r}\gamma^{r}
\end{eqnarray*}
for $i\geq 1,$ where $c=2\gamma/\beta$ and $\gamma$ is as in (\ref{shanxi}). Assuming (\ref{shanxi}), then $(X_1, \cdots, X_k)$ converges weakly to a multivariate normal distribution $N_k(\mu,\Sigma)$ for some $\mu$ and $\Sigma$ given in Theorem 1.5 from \cite{DumitriuGlobal}.
\end{theorem}

Killip \cite{Killip08} obtains the central limit theorem for $\sum_{i=1}^n I(\lambda_i\in (a,b))$, where $a<b$ are two constants. Theorem \ref{hay} is the central limit theorem for homogenous polynomials of $\lambda_i$'s.

Thanks to the recent results of Ram\'{i}rez and Rider \cite{Rider} and Ram\'{i}rez, Rider and Vir\'{a}g \cite{RRV}, we are able to investigate the asymptotic distributions of the smallest and largest eigenvalues for the beta-Jacobi ensembles next.
\noindent Look at the operator
\begin{eqnarray*}
\mathcal{T}_{\beta,a}=-\exp[(a+1)x+\frac{2}{\sqrt{\beta}}b(x)]\frac{d}{dx}\left\{\exp[-ax-\frac{2}{\sqrt{\beta}}b(x)]\right\}
\end{eqnarray*}
where $a>-1$ and $\beta>0$ are constants, and $b(x)$ is a standard Brownian motion on $[0, \infty).$ With probability one, when restricted to the positive half-line with Dirichlet conditions at the origin, $\mathcal{T}_{\beta, a}$ has discrete spectrum comprised of simple eigenvalues $0< \Lambda_0(\beta, a)<\Lambda_0(\beta, a) < \cdots \uparrow \infty$ as stated in Theorem 1 from \cite{Rider}.

For a sequence of pairwise different numbers $a_1, \cdots, a_n$, let $a^{(1)}> a^{(2)} > \cdots > a^{(n)}$ be their order statistic. The following is the limiting distribution of the first $k$ smallest eigenvalues.
\begin{theorem}\lbl{longlive} Let $\lambda_1,\cdots,
\lambda_n$ be random variables with density function $f_{\beta, a_1, a_2}(\lambda)$ as in
 (\ref{beJacobi}). Let $c>0$ be a constant, and $2\beta^{-1}a_1 - n=c$. If $n\to\infty$ and $a_2\to\infty$ such that $n=o(\sqrt{a_2})$, then $(2\beta^{-1}na_2)\cdot (\lambda^{(n)}, \cdots, \lambda^{(n-k+1)})$ converges weakly to $(\Lambda_0(\beta, c), \Lambda_1(\beta, c),\cdots, \Lambda_{k-1}(\beta, c))$.
\end{theorem}
Now look at another random operator
\begin{eqnarray}\lbl{operator}
-\mathcal{H}_{\beta}=\frac{d^2}{dx^2} - x -\frac{2}{\sqrt{\beta}}b_x'
\end{eqnarray}
where $b_x$ is a standard Brownian motion on $[0, +\infty).$ For $\lambda\in \mathbb{R}$ and  function $\psi(x)$ defined on $[0, +\infty)$ with $\psi(0)=0$ and $\int_0^{\infty}\left((\psi')^2 + (1+x)\psi^2\right)\,dx<\infty,$ we say $(\psi, \lambda)$ is an eigenfunction/eigenvalue pair for $-\mathcal{H}_{\beta}$ if $\int_0^{\infty}\psi^2(x)\,dx=1$ and
\begin{eqnarray*}
\psi''(x)=\frac{2}{\sqrt{\beta}}\psi(x)b_x'+(x+ \lambda)\psi(x)
\end{eqnarray*}
holds in the following integration-by-parts sense,
\begin{eqnarray*}
\psi'(x)-\psi'(0)=\frac{2}{\sqrt{\beta}}\psi(x)b_x+\int_0^x-\frac{2}{\sqrt{\beta}}b_y\psi'(y)\,dy + \int_0^x(y +\lambda)\psi(y)\,dy.
\end{eqnarray*}
Theorem 1.1 from \cite{Rider} says that, with probability one, for each $k\geq 1,$ the set of eigenvalues of $-\mathcal{H}_{\beta}$ has well-defined $k$-largest eigenvalues $\Lambda_k.$  Recall (\ref{beJacobi}), set
\begin{eqnarray*}
m_n=\left(\sqrt{n} +\sqrt{2\beta^{-1}a_1}\,\right)^2\ \ \mbox{and}\ \ \sigma_n=\frac{(2\beta^{-1}na_1)^{1/6}}{(\sqrt{n} +\sqrt{2\beta^{-1}a_1}\ )^{4/3}}.
\end{eqnarray*}
Our last result is about the limiting distribution of the first $k$ largest eigenvalues of the beta-Jacobi ensembles.
\begin{theorem}\lbl{work} For each $k\geq 1,$ let $\Lambda_k$ be the $k$-th largest eigenvalue of $-\mathcal{H}_{\beta}$ as in (\ref{operator}). Let $\lambda_1,\cdots,
\lambda_n$ be random variables with joint density function $f_{\beta, a_1, a_2}(\lambda)$ as in
 (\ref{beJacobi}). Assuming (\ref{shanxi}), then $\sigma_{n}\left((2a_2\beta^{-1})\lambda^{(l)}-m_n\right)_{l=1,\cdots,k}$ converges weakly to $(\Lambda_1,\cdots, \Lambda_k)$.
\end{theorem}

\begin{remark} Dumitriu and Koev \cite{Dumitriu2008} derive the exact formulas for the cumulative distribution functions of the largest and smallest eigenvalue of the beta-Jacobi ensembles as in (\ref{beJacobi}) for fixed parameter $\beta, a_1, a_2$ and $n$. Here we complement their work by proving the asymptotic distributions in Theorems \ref{longlive} and \ref{work}.
\end{remark}

\begin{remark}
 In \cite{Jiang07}, Jiang study Theorems \ref{snowf}, \ref{branch}, \ref{hay} and \ref{work} for $\beta=1$ and $2$, which are special cases of the current theorems. The method used in \cite{Jiang07} is the approximation of the entries of Haar-invariant orthogonal or  unitary matrices by independent and identically distributed real or complex Gaussian random variables.
\end{remark}

Now, let us state the methodology used in our proofs. In fact, we employ a different approach than the standard ones in the random matrix theory. Some of the standard tools are the moment method in  \cite{BSJW, Evans, shasha, Jonsson, Wigner}, the Stieltjes transformations in \cite{Lytova, MP67}, the analysis to study the probability density functions of eigenvalues in \cite{adler, BDJ, Johansson, speed, TW, Tracy96, Tracy94}, the large deviation method for obtaining the law of large numbers in \cite{AGZ, Gui, Hiaipetz}, and the application of the free probability theory in \cite{BFOS, BDJiang}. A technique on refined estimates of the smallest eigenvalues is used in \cite{Taocircular}.

 Another way to study random matrices is using the known conclusions, and  connect them with the target ones by approximation. For example, large sample correlation matrices can be approximated by the Wishart matrices \cite{Jiang04}; a large dimensional, Haar-invariant matrix from the classical compact groups  can be approximated by matrices with  independent Gaussian entries \cite{Jiang092, Jiang05}.  The study of local statistics of the Wigner matrices with non-Gaussian entries can be approximated by Gaussian entries \cite{TaoVu}.

In this paper, we approximate the beta-Jacobi ensembles by the beta-Laguerre ensembles through measuring the variation distance between the eigenvalues in the two ensembles (Theorem \ref{main} in Section \ref{yang}). Then the known results for the beta-Laguerre ensemble are used to get those for the beta-Jacobi ensembles.

Now we would like to mention some future problems. Notice that all the theorems above are based on the restriction (\ref{shanxi}). We think it could be relaxed in some situations. One possible way is that, instead of using the {\it uniform} approximation in (\ref{variation}), one can treat case by case for the statistics concerned in the above theorems. For example, to improve Theorem \ref{snowf}, one could  directly evaluate the moments $\sum_{i=1}^n\lambda_i^k$ for $k\geq 1$  by computing the integration with respect to the density function $f_{\beta, a_1, a_2}(\lambda)$  in (\ref{beJacobi}), and then check what restrictions on $a_1, a_2$ and $n$ can make the integral close to the corresponding quantity in the Laguerre case.

Two other problems are discussed in Remarks \ref{gym} and \ref{stool}.

 Finally we give the outline of this paper. In Section \ref{sectionLaguerre}, some known conclusions and some  results on the beta-Laguerre ensembles are reviewed and proved, respectively. They will be used in the proofs for the beta-Jacobi ensembles. In Section \ref{yang}, an approximation theorem for the Jacobi ensembles by the Laguerre ensembles is derived. In Section \ref{shortproof}, we prove the main results stated in this section. In Section \ref{appendix}, some known and useful results are collected for our proofs.

\section{Some Auxiliary Results on $\beta$-Laguerre Ensembles}\lbl{sectionLaguerre}
\setcounter{equation}{0}

Let $\beta>0$ be a constant, $n\geq 2$ be an integer,  $p=1+\frac{\beta}{2}(n-1)$ and
 parameter $a>\frac{\beta}{2}(n-1)$. A $\beta$-Laguerre (Wishart) ensemble is a set of non-negative random variables $(\lambda_1, \cdots, \lambda_n):=\lambda$ with
probability density function
\begin{eqnarray}\lbl{bWishart}
f_{\beta, a}(\lambda)=c_{L}^{\beta, a}\prod_{1\leq i<j\leq
n}|\lambda_i-\lambda_j|^{\beta}\cdot\prod_{i=1}^n\lambda_i^{a-p}\cdot
e^{-\frac{1}{2}\sum_{i=1}^n\lambda_i},
\end{eqnarray}
 where
\begin{eqnarray}\lbl{conWishart}
c_{L}^{\beta,
a}=2^{-na}\prod_{j=1}^{n}\frac{\Gamma(1+\frac{\beta}{2})}
{\Gamma(1+\frac{\beta}{2}j)\Gamma(a-\frac{\beta}{2}(n-j))}.
\end{eqnarray}
One can see \cite{Dumi} for the construction of a matrix to generate eigenvalues with such a distribution. If $\bd{X}=(x_{ij})$ is an $m\times n$ matrix with $m\geq n$, where $x_{ij}$'s are independent and identically distributed random variables with the standard real normal ($\beta=1$), complex normal ($\beta=2$) or quaternion normal ($\beta=4$) distribution,
then $f_{\beta}(\lambda)$ is the density function of the eigenvalues
$\lambda=(\lambda_1, \cdots, \lambda_n)$ of $\bd{X}^*\bd{X}$ with
$a=\frac{\beta}{2}m$ for $\beta=1,2,$ or $4$. See
\cite{Edelman1989, James1964, Muirhead1982} for the cases $\beta=1$ and $2$,
and \cite{Macdonald1998} for $\beta=4$, or (4.5) and (4.6) from
\cite{Edelman05}.

It is easy to see from Theorem A.\ref{haci} in Appendix that the following is true.

\begin{lemma}\lbl{flu} Let  $\lambda=(\lambda_1,\cdots,
\lambda_n)$ be random variables with the density function as in
(\ref{bWishart}). If $n\to \infty,\ a\to \infty$ and $n\beta/(2a)\to \gamma\leq 1$, then
\begin{eqnarray*}
& (i) & \frac{1}{n^{i+1}}\sum_{j=1}^n\lambda_j^i \ \mbox{converges to}\
\left(\frac{\beta}{\gamma}\right)^i\sum_{r=0}^{i-1}\frac{1}{r+1}\binom{i}{r}\binom{i-1}{r}\gamma^{r}\ \ \mbox{in probability};\\
& (ii)  & \frac{1}{n^{i}}\sum_{j=1}^n\lambda_j^i -
\left(\frac{\beta}{\gamma}\right)^i n\sum_{r=0}^{i-1}\frac{1}{r+1}\binom{i}{r}\binom{i-1}{r}\gamma^{r}\ \mbox{converges to}\  N(\mu_i, \sigma_i^2)
\end{eqnarray*}
in distribution for any integer $i\geq 1,$  where $\mu_i$ and $\sigma_i^2$ are constants depending on $\gamma, \beta$ and $i$ only.
\end{lemma}

A large deviation result in \cite{Hiaipetz, Hiai} implies the following ``law of large numbers".
\begin{lemma}\lbl{xiaochild} Let  $\lambda_1,\cdots,
\lambda_n$ be random variables with the density  as in
(\ref{bWishart}). Assume $n\beta/(2a)\to \gamma\in (0, 1],$ and let
$\gamma_{min}=(\sqrt{\gamma}-1)^2$ and
$\gamma_{max}=(\sqrt{\gamma}+1)^2$. Let $\mu_n$ be the empirical distribution of $Y_i:=\lambda_i\gamma/(n\beta)$ for $i=1,2, \cdots,n.$ Then $\mu_n$ converges weakly to the distribution $\mu_{\infty}$ with density $f_{\gamma}(x)$ as in (\ref{mooncake}) almost surely. Moreover, $\liminf_{n\to\infty}Y_{max}(n)\geq \gamma_{max}\ a.s.$ and $\limsup_{n\to\infty}Y_{min}(n)\leq \gamma_{min}\ a.s.,$ where $Y_{max}(n)=\max\{Y_1, \cdots, Y_n\}$ and $Y_{min}(n)=\min\{Y_1, \cdots, Y_n\}$.
\end{lemma}
The first part of the above lemma is also obtained in \cite{Dumitriu1}.

\noindent\textbf{Proof of Lemma \ref{xiaochild}}. From (\ref{bWishart}), it is obvious that the joint density function of $\{Y_i,\ 1\leq i \leq n\}$ is
\begin{eqnarray*}
g(y_1, \cdots, y_n)=Const\cdot \prod_{1\leq i < j\leq n}|y_i - y_j|^{\beta}\, \prod_{i=1}^ny_i^{a-p}\, e^{-(n\beta/2\gamma)\sum_{i=1}^ny_i}
\end{eqnarray*}
for $y_1>0,\cdots, y_n>0.$ Take $\gamma(n)=a-p,\, p(n)=n, Q(t)=\beta t/(2\gamma)$ in Theorem A.\ref{casino}. Then $p(n)/n= 1$ and $\gamma(n)/n=(a-p)/n\to \beta (\gamma^{-1}-1)/2$ as $n\to\infty$ since $p=1+\beta(n-1)/2.$ According to Theorem  A.\ref{casino}, $\mu_n$ satisfies the large deviation principle with speed $\{1/n^2;\, n\geq 1\}$ and rate function
\begin{eqnarray*}
I(\nu) &= & -\frac{\beta}{2}\iint\log |x-y|d\nu(x)d\nu(y) + \int\left(\frac{\beta x}{2\gamma}-\frac{\beta(\gamma^{-1}-1)}{2}\log x\right)d \nu(x) + B\\
& = & B + \frac{\beta}{\gamma^2}\left(-\frac{\gamma^2}{2}\iint\log |x-y|\, d\nu(x)\,d\nu(y) + \frac{\gamma}{2}\int\left(x-(1-\gamma)\log x\right)d\nu(x)\right)
\end{eqnarray*}
where $B$ is a finite constant. By Theorem A.\ref{ok}, the unique measure $\mu_{\infty}$  on $[0, \infty)$ to achieve the minimum of $I(\nu)$ over all probability measures on $[0, \infty)$ is the Marchenko-Pastur law with density function $f_{\gamma}(x)$ as in (\ref{mooncake}).

For $\epsilon>0,$ let $F=\{\nu;\, d(\nu, \mu_{\infty})\geq \epsilon\}$, where $d$ is as in (\ref{du}). Then, $F$ is a closed set in the weak topology. By the large deviation upper bound, there exists a constant $C>0$ such that $P(d(\mu_n, \mu_{\infty}) >\epsilon)\leq e^{-n^2C}$ as $n$ is large enough. By the Borel-Cantelli lemma, $d(\mu_n, \mu_{\infty})\to 0\ a.s.$ as $n \to\infty.$ The first part of the conclusion is proved.

For any integer $k\geq 1,$ it is easy to see
\begin{eqnarray*}
Y_{max}(n) \geq \left(\int_{0}^{\infty}y^k\,d\mu_n(y)\right)^{1/k}\ \ \mbox{and}\ \ \frac{1}{Y_{min}(n)} \geq \left(\int_{0}^{\infty}y^{-k}\,d\mu_n(y)\right)^{1/k}.
\end{eqnarray*}
Since $\mu_n$ converges weakly to $\mu_{\infty}$ almost surely, by the Fatou lemma,
\begin{eqnarray*}
& & \liminf_{n\to\infty}Y_{max}(n)\geq \left(\int_{\gamma_{min}}^{\gamma_{max}}y^{k}f_{\gamma}(y)\,dy\right)^{1/k}\ a.s.\ \ \mbox{and}\\
& & \liminf_{n\to\infty}\frac{1}{Y_{min}(n)}\geq \left(\int_{\gamma_{min}}^{\gamma_{max}}y^{-k}f_{\gamma}(y)\,dy\right)^{1/k}\ a.s.
\end{eqnarray*}
for any integer $k\geq 1.$ Letting $k\to \infty,$ we have
\begin{eqnarray*}
\liminf_{n\to\infty}Y_{max}(n)\geq \gamma_{max}\ a.s.\ \ \mbox{and}\ \ \limsup_{n\to\infty}Y_{min}(n)\leq \gamma_{min}\ a.s.\ \ \ \ \ \ \ \blacksquare
\end{eqnarray*}
In what follows, the notation $\chi^2(s)$ stands for the $\chi^2$ distribution with degrees of freedom $s$; $\chi(s)$ denotes a positive random variables with $(\chi(s))^2$ following the $\chi^2(s)$ distribution.

\begin{lemma}\lbl{chimono} Let $X$ have the Gamma distribution with density function
\begin{eqnarray*}
f(x|\alpha,\theta)=\begin{cases}
\frac{x^{\alpha-1}e^{-x/\theta}}{\Gamma(\alpha)\, \theta^{\alpha}}, & \text{if $x\geq 0$;}\\
0, & \text{otherwise,}
\end{cases}
\end{eqnarray*}
 where $\alpha>0$ and $\theta>0$ are constants. Given $b>0$ and $\theta>0$, set $g(\alpha):=P(X\geq b)$. Then $g(\alpha)$ is increasing over $(0, \infty).$ In particular, $P(\chi^2(s) \geq b)\leq P(\chi^2(t) \geq b)$ for any $0<s<t$ and $b>0.$
\end{lemma}
\noindent \textbf{Proof}. Since $\chi^2(s)$ has probability density function $f(x|\alpha,\theta)$ with $\alpha=p/2$ and $\theta=2$, we only need to prove the first part of the theorem.

Obviously, $X/\theta$ has density function $f(x|\alpha, 1)$, without loss of the generality, we prove the conclusion by assuming $\theta=1.$ First, noticing $\Gamma(\alpha)=\int_0^\infty x^{\alpha-1}e^{-x}\,dx$ for $\alpha>0,$ then
\begin{eqnarray*}
g_1(\alpha):=g(\alpha + 1)=\frac{\int_{b}^{\infty}x^{\alpha}e^{-x}\,dx}{\int_0^\infty x^{\alpha}e^{-x}\,dx}
\end{eqnarray*}
for $\alpha>-1$. Use the fact that $\frac{d\, x^{\alpha}}{d\alpha}=x^{\alpha}\log x$ for any $x>0$ and $\alpha \in \mathbb{R}$ to  obtain that
\begin{eqnarray*}
\frac{d g_1(\alpha)}{d\alpha}&= & \frac{\int_{b}^{\infty}x^{\alpha}(\log x)e^{-x}\,dx\cdot \int_0^\infty x^{\alpha}e^{-x}\,dx-  \int_{b}^{\infty}x^{\alpha}e^{-x}\,dx\cdot \int_0^\infty x^{\alpha}(\log x)e^{-x}\,dx }{(\int_0^\infty x^{\alpha}e^{-x}\,dx)^2}\\
 & = & \frac{\int_0^\infty x^{\alpha}e^{-x}\,dx\cdot \int_b^\infty x^{\alpha}e^{-x}\,dx}{\Gamma(\alpha + 1)^2}(h(b)-h(0)),
\end{eqnarray*}
where $h(b):=\int_b^{\infty}x^{\alpha}(\log x) e^{-x}\,dx/\int_b^\infty x^{\alpha}e^{-x}\,dx$ for $b\geq 0.$ Fix $\alpha>-1$. To show $\frac{d g_1(\alpha)}{d\alpha}>0$ for all $\alpha>-1$, it is enough to show that $h(b)$ is strictly increasing on $b\in [0,\infty).$ In fact,
\begin{eqnarray*}
& & \frac{d\, h(b)}{db}\\
& = &\frac{1}{(\int_b^\infty x^{\alpha}e^{-x}\,dx)^2}\left(-b^{\alpha}e^{-b}(\log b)\int_b^\infty x^{\alpha}e^{-x}\,dx +  b^{\alpha}e^{-b}\int_b^\infty x^{\alpha}(\log x)e^{-x}\,dx\right)\\
&>& 0
\end{eqnarray*}
since $\int_b^\infty x^{\alpha}(\log x)e^{-x}\,dx> (\log b)\int_b^\infty x^{\alpha}e^{-x}\,dx$ for all $b>0.$\ \ \ \ \ \ $\blacksquare$\\


\begin{lemma}\lbl{glue} Let $\{a_n;\ n\geq 1\}$ and $\beta$ be positive constants with $2a_n>(n-1)\beta$ for all  $n\geq 1$. Set $\mathbb{T}_n=\left\{(n-i+1)\beta,\ 2a_n-(i-1)\beta;\ 1\leq i \leq n\right\}$ for $n\geq 1.$ Let $\{\chi(t),\ t\in \mathbb{T}_n;\ n\geq 1\}$ be a set of random variables defined on the same probability space.
Define
\begin{eqnarray}
\delta_n= \max_{t\in\mathbb{T}_n}|\chi(t) - \sqrt{t} |\ \ \ \mbox{and}\ \ \ \rho_n= \max_{t\in\mathbb{T}_n}|\chi^2(t) - t |
\end{eqnarray}
for $n \geq 1$. If $n\beta/(2a_n)\to \gamma\in (0, 1]$, then $\delta_n/\sqrt{n}\to 0\ a.s.$ and $\rho_n/n\to 0\ a.s.$ as $n\to\infty$.
\end{lemma}

\noindent\textbf{Proof}. First, noticing $|\mathbb{T}_n|\leq 2n$ and $\max\{t;\, t\in \mathbb{T}_n\}\leq n\beta + 2a_n$ for all $n\geq 1$. Since $n\beta/(2a_n)\to \gamma$ as $n\to \infty,$ there exists a constant $C\geq 1$ such that $\max\{t;\, t\in \mathbb{T}_n\}\leq Cn$ for all $n\geq 1.$
Then, for any $\epsilon>0$,
\begin{eqnarray}\lbl{aaa}
P\left(\delta_n \geq \sqrt{n}\epsilon\right)
& \leq & 2n \cdot\max_{t\in \mathbb{T}_n}P\left(|\chi(t) - \sqrt{t} |\geq \sqrt{n}\epsilon\right)\nonumber\\
& \leq & 2n \cdot\max_{0\leq t \leq Cn}P\left(|\chi(t) - \sqrt{t} |\geq \sqrt{n}\epsilon\right)
\end{eqnarray}
for all $n\geq 1$.  Now, taking $p_n=[n^{1/3}]$ for $n\geq 1,$ by Lemma \ref{chimono},
\begin{eqnarray*}
\max_{0\leq t \leq p_n}P\left(|\chi(t) - \sqrt{t} |\geq \sqrt{n}\epsilon\right)
& \leq & \max_{0\leq t \leq p_n}P\left(\chi(t) \geq \sqrt{n}\epsilon-\sqrt{p_n}\right)\\
 & \leq & P\left(\chi(p_n) \geq \sqrt{n}\epsilon-\sqrt{p_n}\right)\leq P\left(|\chi(p_n) -\sqrt{p_n}| \geq \frac{\sqrt{n}\epsilon}{2}\right)
\end{eqnarray*}
as $n$ is large enough since $p_n=o(\sqrt{n})$ as $n\to\infty.$  
This and (\ref{aaa}) imply that, given $\epsilon>0,$
\begin{eqnarray}\lbl{juice}
P\left(\delta_n \geq \sqrt{n}\epsilon\right) \leq 2n\cdot\max_{p_n\leq t \leq Cn}P\left(|\chi(t) - \sqrt{t}|\geq \frac{\sqrt{n}\epsilon}{2}\right)
\end{eqnarray}
as $n$ is sufficiently large. Now, the last probability is equal to
\begin{eqnarray*}
& & P\left(\chi(t) \geq \sqrt{t} + \frac{\sqrt{n}\epsilon}{2}\right) + P\left(\chi(t) \leq \sqrt{t} - \frac{\sqrt{n}\epsilon}{2}\right)\\
& \leq & P\left(\chi([t]+1) \geq \sqrt{t} + \frac{\sqrt{n}\epsilon}{2}\right) + P\left(\chi([t]) \leq \sqrt{t} - \frac{\sqrt{n}\epsilon}{2}\right)\ \ \ \ (\mbox{by Lemma \ref{chimono}})\\
& \leq & P\left(\chi([t]+1) - \sqrt{[t]+1} \geq -1 + \frac{\sqrt{n}\epsilon}{2}\right) + P\left(\chi([t])- \sqrt{[t]} \leq 1 - \frac{\sqrt{n}\epsilon}{2}\right)\\
& \leq & 2\, P\max_{[t]\leq k \leq [t]+1}P(|\chi(k)-\sqrt{k}|\geq \frac{\sqrt{n}\epsilon}{3})
\end{eqnarray*}
as $n$ is sufficiently large, where the inequality $\sqrt{1+[t]}-1 \leq \sqrt{t} \leq 1+\sqrt{[t]}$ for all $t\geq 1$ is used in the last step. Since  $|\chi(k)-\sqrt{k}|\leq |(\chi(k))^2-k|/\sqrt{k}$ for any $k\geq 1,$ by (\ref{juice}),
\begin{eqnarray}\lbl{bao}
P\left(\delta_n \geq \sqrt{n}\epsilon\right) \leq 4n\cdot\max_{p_n\leq k \leq 2Cn}P\left(\frac{|(\chi(k))^2 - k|}{\sqrt{k}}\geq \frac{\sqrt{n}\epsilon}{3}\right)
\end{eqnarray}
as $n$ is sufficiently large. Lemma 2.4 from \cite{Jiang09} says that
\begin{eqnarray*}
P\left(\frac{|\sum_{i=1}^n(\xi_i^2-1)|}{\sqrt{n}}\geq
c\right)\leq 2e^{-c^2/6}
\end{eqnarray*}
for any $n\geq 1$ and $c \in (0,\, \sqrt{n}/2)$, where  $\xi_1, \xi_2, \cdots, \xi_n$ are i.i.d.
random variables with $\xi_1\sim N(0,1)$. Since $\sqrt{n}\epsilon/3\geq \sqrt{p_n}/3$ as $n$ is large enough, the above inequality implies that
\begin{eqnarray*}
\max_{p_n\leq k \leq 2Cn}P\left(\frac{|(\chi(k))^2 - k|}{\sqrt{k}}\geq \frac{\sqrt{n}\epsilon}{3}\right)
&\leq & \max_{p_n\leq k \leq 2Cn}P\left(\frac{|(\chi(k))^2 - k|}{\sqrt{k}}\geq \frac{\sqrt{p_n}}{3}\right)\\
& \leq & 2\, e^{-n^{1/3}/55}
\end{eqnarray*}
as $n$ is sufficiently large. This and (\ref{bao}) tell us that, for any $\epsilon >0$, $P\left(\delta_n \geq \sqrt{n}\epsilon\right)< e^{-n^{1/3}/56}$ as $n$ is large enough. Therefore, by the Borel-Cantelli lemma, $\delta_n/\sqrt{n} \to 0\ \ a.s. $ as $n\to\infty$.

Finally, since $a^2-b^2=(a-b)^2 + 2b(a-b)$ for any $a, b\in \mathbb{R},$ recalling the constant $C\geq 1$ from (\ref{aaa}), we have
\begin{eqnarray*}
\frac{\rho_n}{n}
& \leq & \max_{t\in \mathbb{T}_n}\frac{|\chi(t)-\sqrt{t}|^2}{n} + \max_{t\in \mathbb{T}_n}\left\{\frac{2\sqrt{t}}{\sqrt{n}}\cdot\frac{|\chi(t)-\sqrt{t}|}{\sqrt{n}}\right\}\\
& \leq & \left(\frac{\delta_n}{\sqrt{n}}\right)^2 + (2C)\cdot\frac{\delta_n}{\sqrt{n}}\to 0\ \ a.s.
\end{eqnarray*}
as $n\to\infty$. \ \ \ \ \ \ \ $\blacksquare$

By following the proof in \cite{Jack}, we have a result below on the $\beta$-Laguerre ensembles. It is also reported in Theorem 10.2.2 from \cite{Dumitriu1} without proof.

\begin{lemma}\lbl{repair} Let  $\lambda=(\lambda_1,\cdots,
\lambda_n)$ be random variables with the density function as in
(\ref{bWishart}). Set $\lambda_{max}(n)=\max\{\lambda_1,\cdots, \lambda_n\}$ and $\lambda_{min}=\min\{\lambda_1,\cdots, \lambda_n\}$. If $n\to \infty,\ a\to \infty$ and $n\beta/(2a)\to \gamma\in (0, 1]$, then
\begin{eqnarray*}
\frac{\lambda_{max}(n)}{n} \to \beta (1+ \sqrt{\gamma^{-1}})^2\ a.s.\ \ \mbox{and}\ \ \  \frac{\lambda_{min}(n)}{n} \to \beta (1-\sqrt{\gamma^{-1}})^2\ \ a.s.
\end{eqnarray*}
 as $n\to \infty.$
\end{lemma}
\noindent \textbf{Proof}. By Lemma \ref{xiaochild}, it is enough to show
\begin{eqnarray}
& & \limsup_{n\to\infty} \frac{\lambda_{max}(n)}{n} \leq \beta (1+\sqrt{\gamma^{-1}})^2\ \ a.s.\lbl{ming}\ \ \mbox{and}\\
& & \liminf_{n\to\infty} \frac{\lambda_{min}(n)}{n} \geq \beta (1-\sqrt{\gamma^{-1}})^2\ \ a.s.\lbl{yue}
\end{eqnarray}
 From Theorem 3.4 in \cite{Dumi}, we know the eigenvalues of $B_{\beta}B_{\beta}^T$ have the same probability distribution as in (\ref{bWishart}), where
\begin{eqnarray*}
B_{\beta}=
\begin{pmatrix}
\chi_{2a} & & & & \\
\chi_{\beta(n-1)} & \chi_{2a-\beta} & & \\
& \ddots & \ddots & \\
& & \chi_{\beta} & \chi_{2a-\beta (n-1)}
\end{pmatrix}
_{n\times n}
\end{eqnarray*}
where all of the $2n-1$ entries in the matrix are non-negative independent random variables with  $(\chi_{s})^2$ following  the $\chi^2$ distribution with degree of freedom $s$.
It is easy to see that the first and last rows of $B_{\beta}B_{\beta}^T$ are respectively are
\begin{eqnarray*}
(\chi_{2a}^2, \,\chi_{2a}\chi_{\beta(n-1)},0,\cdots,0, 0)\ \ \mbox{and}\ \ (0,0,0,\cdots,\chi_{\beta} \chi_{2a-\beta (n-2)} ,\, \chi_{\beta}^2 + \chi_{2a-\beta (n-1)}^2).
\end{eqnarray*}
For $i=2,\cdots, n-1$ and $n\geq 3,$ the $i$-th row is
\begin{eqnarray*}
(0, \cdots,0,\chi_{\beta(n-i+1)}\chi_{2a-(i-2)\beta},\underbrace{\chi_{\beta(n-i+1)}^2 + \chi_{2a-(i-1)\beta}^2}_{ith \, position}, \chi_{\beta(n-i)}\chi_{2a-(i-1)\beta},0,\cdots,0).
\end{eqnarray*}
Recall the Ger\^{s}gorin theorem: each eigenvalue of an $n\times n$ matrix $A=(a_{ij})$ lies in at least one of the disks $\{z\in \mathbb{C};\ |z-a_{jj}|\leq \sum_{i\ne j}|a_{ij}|\}$ for $j=1,2,\cdots,n.$    Then
\begin{eqnarray}\lbl{kang}
 \lambda_{max}(n) \leq \max\{U^+, V^+, W_i^+;\ 2\leq i \leq n-1\}
\end{eqnarray}
where
\begin{eqnarray}
& & U^+=\chi_{2a}^2 + \chi_{2a}\chi_{\beta(n-1)},\ \ V^+=\chi_{\beta}^2 + \chi_{2a-\beta (n-1)}^2 + \chi_{\beta} \chi_{2a-\beta (n-2)}\ \ \ \mbox{and}\nonumber\\
& & W_i^+=\chi_{\beta(n-i+1)}^2 + \chi_{2a-(i-1)\beta}^2  + \chi_{\beta(n-i+1)}\chi_{2a-(i-2)\beta} + \chi_{\beta(n-i)}\chi_{2a-(i-1)\beta}\lbl{sud}
\end{eqnarray}
for $2\leq i \leq n-1$,
and
\begin{eqnarray}\lbl{eds}
\lambda_{min}(n) \geq \min\{U^-, V^-, W_i^-;\ 2\leq i \leq n-1\}
\end{eqnarray}
where
\begin{eqnarray*}
& & U^-=\chi_{2a}^2 -\chi_{2a}\chi_{\beta(n-1)},\ V^-=\chi_{\beta}^2 + \chi_{2a-\beta (n-1)}^2 - \chi_{\beta} \chi_{2a-\beta (n-2)}\ \ \mbox{and}\\
& & W_i^-=\chi_{\beta(n-i+1)}^2 + \chi_{2a-(i-1)\beta}^2  - \chi_{\beta(n-i+1)}\chi_{2a-(i-2)\beta} - \chi_{\beta(n-i)}\chi_{2a-(i-1)\beta}.
\end{eqnarray*}
Set
\begin{eqnarray*}
Q_{n,i}^+= \beta(n-i+1) + 2a-(i-1)\beta &+ & \sqrt{\beta(n-i+1)(2a-(i-2)\beta)} \\
& + & \sqrt{\beta(n-i)(2a-(i-1)\beta)}
\end{eqnarray*}
for $2\leq i \leq n-1$ and $n\geq 3.$ Reviewing the notation $\delta_n$ and $\mathbb{T}_n$ in Lemma \ref{glue}, we have  that
\begin{eqnarray}
& & |\chi_{\beta(n-i+1)}\chi_{2a-(i-2)\beta} - \sqrt{\beta(n-i+1)(2a-(i-2)\beta)}\, |\nonumber\\
& = & |\,\left(\chi_{\beta(n-i+1)} -\sqrt{\beta(n-i+1)}\right)\left(\chi_{2a-(i-2)\beta} - \sqrt{2a-(i-2)\beta}\right)  \nonumber\\
& & + \sqrt{\beta(n-i+1)}\cdot \left(\chi_{2a-(i-2)\beta} - \sqrt{2a-(i-2)\beta}\right)\nonumber \\
& & + \sqrt{2a-(i-2)\beta}\left(\chi_{\beta(n-i+1)} -\sqrt{\beta(n-i+1)}\right)\,|\nonumber\\
& \leq & \delta_n^2 + (\sqrt{\beta n} + \sqrt{2a})\delta_n\lbl{def}
\end{eqnarray}
uniformly for all $2\leq i \leq n-1$ and $n\geq 3.$ Moreover,
\begin{eqnarray}\lbl{vic}
|\chi_{\beta(n-i+1)}^2-\beta(n-i+1)|\leq \rho_n\ \ \mbox{and}\ \ |\chi_{2a-(i-1)\beta}^2- (2a-(i-1)\beta)|\leq \rho_n
\end{eqnarray}
for all $1\leq i \leq n$ and $n\geq 1.$ Thus,
\begin{eqnarray}\lbl{class}
|\,\max_{2\leq i \leq n-1}\{W_i^+\} - \max_{2\leq i \leq n-1}\{Q_{n,i}^+\}| \leq \max_{2\leq i \leq n-1}|W_i^+ - Q_{n,i}^+| \leq K(\rho_n + \delta_n^2 + \delta_n\sqrt{n})
\end{eqnarray}
as $n$ is sufficiently large, $K$ is a constant not depending on $i$ or $n$. Evidently, $Q_{n,i}^+\leq \beta n + 2a+ 2 \sqrt{2a\beta n}=(\sqrt{\beta n}+ \sqrt{2a}\,)^2$ uniformly for all $2\leq i \leq n-1$ and $n\geq 3.$ By Lemma \ref{glue} and the condition that $n\beta/(2a)\to \gamma$, we have
\begin{eqnarray}\lbl{visitor1}
\limsup_{n\to\infty} \frac{\max_{2\leq i \leq n-1}\{W_i^+\}}{n} \leq \beta (1+\sqrt{\gamma^{-1}})^2\ \ a.s.
\end{eqnarray}
By the same but easier argument, the above is also true when $\max_{2\leq i \leq n-1}\{W_i^+\}$ is replaced by $U^+$ and $V^+$, respectively. Therefore, (\ref{ming}) is concluded from (\ref{kang}) and (\ref{visitor1}).

Now we prove (\ref{yue}). Write $(2a-(i-1))\beta =-\beta + (2a-(i-2))\beta $. It is easy to see
\begin{eqnarray*}
Q_{n,i}^- :
&= & \beta(n-i+1) + (2a-(i-1))\beta -  \sqrt{\beta(n-i+1)(2a-(i-2)\beta)} \\
& & \ \ \ \ \ \ \ \ \ \ \ \ \ \ \ \ \ \ \ \ \ \ \ \ \ \ \ \ \ \ \ \ \ \ \ \ \ \ \ \ - \sqrt{\beta(n-i)(2a-(i-1)\beta)}\\
& \geq & -\beta + \beta(n-i+1) + (2a-(i-2)\beta) -  2\sqrt{\beta(n-i+1)(2a-(i-2)\beta)}\\
& = & -\beta + (\sqrt{\beta(n-i+1)} - \sqrt{2a-(i-2)\beta})^2
\end{eqnarray*}
for $2\leq i \leq n-1$ and $n\geq 3.$ Use the equality $\sqrt{x}-\sqrt{y}=(x-y)/(\sqrt{x}+ \sqrt{y})$ to obtain
\begin{eqnarray*}
(\sqrt{\beta(n-i+1)} - \sqrt{2a-(i-2)\beta})^2 = \frac{(\beta (n-1) -2a)^2}{(\sqrt{\beta(n-i+1)} + \sqrt{2a-(i-2)\beta}\,)^2 },
\end{eqnarray*}
which is increasing in $i\in \{2,3,\cdots,n-1\}$. Since $a> \beta(n-1)/2$, we get
\begin{eqnarray*}
\min_{2\leq i\leq n-1}(\sqrt{\beta(n-i+1)} - \sqrt{2a-(i-2)\beta})^2 = \frac{(\beta (n-1) -2a)^2}{(\sqrt{\beta(n-1)} + \sqrt{2a}\,)^2 }.
\end{eqnarray*}
for $n\geq 3.$ Combining all the above steps, and using the condition that $n\beta/(2a)\to \gamma$, we arrive at
\begin{eqnarray}\lbl{laid}
\liminf_{n\to\infty}\frac{\min_{2\leq i \leq n-1}\{Q_{n,i}^-\}}{n}\geq \frac{\beta(1-\gamma^{-1})^2}{(1+\sqrt{\gamma^{-1}})^2}=\beta (1-\sqrt{\gamma^{-1}})^2\ \ a.s.
\end{eqnarray}
By  (\ref{def}) and (\ref{vic}) and the same argument as in (\ref{class}), we obtain
\begin{eqnarray*}
|\,\min_{2\leq i \leq n-1}\{W_i^-\} - \min_{2\leq i \leq n-1}\{Q_{n,i}^-\}| \leq \max_{2\leq i \leq n-1}|W_i^- - Q_{n,i}^-| \leq K\left(\rho_n + \delta_n^2 + \sqrt{n}\delta_n\right)
\end{eqnarray*}
as $n$ is sufficiently large. This together with (\ref{laid}) and Lemma \ref{glue} concludes
\begin{eqnarray*}
\liminf_{n\to\infty}\frac{\min_{2\leq i\leq n-1}\{W_i^-\}}{n}\geq \beta (1-\sqrt{\gamma^{-1}})^2\ \ a.s.
\end{eqnarray*}
Similarly,
\begin{eqnarray*}
& & \lim_{n\to\infty}\frac{U^-}{n}=\lim_{n\to\infty}\frac{2a-\sqrt{2a\beta(n-1)}}{n}=\beta\gamma^{-1} -\beta{\sqrt{\gamma^{-1}}}\ a.s. \\
& & \lim_{n\to\infty}\frac{V^-}{n}=\lim_{n\to\infty}\frac{2a-\beta(n-1)}{n}=\beta\gamma^{-1} -\beta\ a.s.
\end{eqnarray*}
Since $\gamma \in (0,1],$ it is obvious that $\min\{\beta\gamma^{-1}-\beta \sqrt{\gamma^{-1}},\, \beta\gamma^{-1} -\beta\}\geq \beta (1-\sqrt{\gamma^{-1}})^2.$ By (\ref{eds}), the above three assertions imply (\ref{yue}). \ \ \ \ \ \ $\blacksquare$

\section{An Approximation Result}\lbl{yang}
\setcounter{equation}{0}

Let $\mu$ and $\nu$ be probability measures on $(\mathbb{R}^m, \mathcal{B}),$ where $m\geq 1$ and $\mathcal{B}$ is the Borel $\sigma$- algebra on $\mathbb{R}^m$. The variation distance $\|\mu - \nu\|$ is defined by
\begin{eqnarray}\lbl{variation}
\|\mu- \nu\|= 2\cdot \sup_{A \in \mathcal{B}}|\mu(A) - \nu (A)|= \int_{\mathbb{R}^m}|f(x)-g(x)|\,dx_1\cdots dx_m
\end{eqnarray}
if $\mu$ and $\nu$ have density functions $f(x)$ and $g(x)$ with respect to the Lesbegue measure. For a random variable $Z$, we use $\mathcal{L}(Z)$ to denote its probability distribution. The following is the key tool to prove the results stated in Introduction.

\begin{theorem}\lbl{main} Let $\mu=(\mu_1, \cdots, \mu_n)$ and $\lambda=(\lambda_1, \cdots,
\lambda_n)$ be random variables with density $f_{\beta,
a_1}(\mu)$ as in (\ref{bWishart})(taking $a=a_1$) and $f_{\beta, a_1, a_2}(\lambda)$ as in
 (\ref{beJacobi}). Assuming (\ref{shanxi}), then
$\|\ml{L}(2a_2\lambda)-\ml{L}(\mu)\| \to 0$.
\end{theorem}

 \begin{remark}\lbl{gym} The condition (\ref{shanxi}) is actually sharp for $\|\ml{L}(2a_2\lambda)-\ml{L}(\mu)\| \to 0$ for $\beta=1$ and $2.$ In fact, let $U$ be an $N\times N$ random unitary matrix with real ($\beta=1$) and complex ($\beta=2$), chosen with Haar measure. Decompose
\begin{eqnarray*}
U=\begin{pmatrix}
A_{n_1,n_2} & C_{n_1\times (N-n_2)}\\
B_{(N-n_1)\times n_2} & D_{(N-n_1)\times (N-n_2)}
\end{pmatrix}
.
\end{eqnarray*}
Assume $n_1\geq n_2$ and $n_1+n_2\leq N.$ Set $Y=A_{n_1,n_2}^*A_{n_1,n_2}$, we see from (3.15) in \cite{Conductance}  that the eigenvalues $\lambda_1,\cdots,\lambda_{n_2}$ of $Y$ have density function
\begin{eqnarray*}
Const\cdot \prod_{1\leq
i<j\leq
n}|\lambda_i-\lambda_j|^{\beta}\cdot\prod_{i=1}^n\lambda_i^{a_1-p}(1-\lambda_i)^{a_2-p},
\end{eqnarray*}
which belongs to the beta-Jacobi ensemble (\ref{beJacobi}) with
\begin{eqnarray*}
a_1=\frac{\beta}{2}n_1,\ \  a_2=\frac{\beta}{2}(N-n_1),\  \ n=n_2\  \ \mbox{and}\  \ p= 1+ \frac{\beta}{2}(n-1).
\end{eqnarray*}
It is shown in \cite{Jiang06} for $\beta=1$ and \cite{Jiang07} for $\beta=2$ that the variation distance $d'$ between $\sqrt{N}A_{n_1,n_2}$ and $X_n$ goes to zero if $n_1=o(\sqrt{N})$ and $n_2=o(\sqrt{N})$, where $X_n$ is the joint distribution of $n_1n_2$ independent and identically distributed real Gaussian ($\beta=1$) or complex Gaussian ($\beta=2$) random variables. The orders $n_1=o(\sqrt{N})$ and $n_2=o(\sqrt{N})$, which correspond to that $a_1=o(\sqrt{a_2})$ and $n_1=o(\sqrt{a_2})$ in (\ref{shanxi}), are also proved to be sharp for both cases in \cite{Jiang07, Jiang06}. A further analysis shows that $d'\to 0$ if and only if $\|\ml{L}(2a_2\lambda)-\ml{L}(\mu)\| \to 0$. This tells us that the orders in Theorem \ref{main} are sharp for $\beta=1$ and $2.$ However, it is not known whether the same remains true for other $\beta>0$.
\end{remark}

\begin{remark}\lbl{stool} The condition ``$n\beta/(2a_1)\to \gamma$" in (\ref{shanxi}) is required in Theorem \ref{main}. In the same contexts as in \cite{Jiang06} for $\beta=1$ and \cite{Jiang07} for $\beta=2$, the condition is not imposed. Although the condition is harmless in proving the main results stated in Introduction, it would be interesting to see if it can be removed for other $\beta>0.$
\end{remark}

Now we start to prove the above theorem by developing several lemmas.
\begin{lemma}\lbl{botany} Let $n\geq 2.$ Let $\mu=(\mu_1, \cdots, \mu_n)$ and $\lambda=(\lambda_1, \cdots,
\lambda_n)$ be random variables with density functions $f_{\beta,
a_1}(\mu)$ as in (\ref{bWishart})
(taking $a=a_1$) and $f_{\beta, a_1, a_2}(\lambda)$  in (\ref{beJacobi}), respectively. Then
\begin{eqnarray*}
\|\ml{L}(2a_2\lambda)-\ml{L}(\mu)\|=E|K_{n}\cdot
L_n(\mu)-1|
\end{eqnarray*}
and $E(K_{n}\cdot
L_n(\mu))=1$,  where
\begin{eqnarray}
& & K_{n}=a_2^{-na_1}\cdot
\prod_{i=0}^{n-1}\frac{\Gamma(a_1+a_2-\frac{\beta}{2}i)}
{\Gamma(a_2-\frac{\beta}{2}i)}\ \ \mbox{and}\lbl{shining1}\\
& & L_n(\mu)=e^{(1/2)\sum_{i=1}^n{\mu_i}}\cdot
\prod_{i=1}^n\left(1-\frac{\mu_i}{2a_2}\right)^{a_2-p}
 \cdot I(\max_{1\leq
i \leq n}\mu_i \leq 2a_2).\lbl{shining2}
\end{eqnarray}
\end{lemma}

\noindent\textbf{Proof.} It is enough to show
\begin{eqnarray}\lbl{celcius}
\|\ml{L}(2a_2\lambda)-\ml{L}(\mu)\|=\int_{[0,\infty)^n} |K_{n}\cdot
L_n(\mu)-1|\cdot f_{\beta, a_1}(\mu)\,d\mu.
\end{eqnarray}
First, since $p=\frac{\beta}{2}(n-1) + 1$, we have $n(n-1)\beta/2+n(a_1-p)+n=na_1.$ It is easy to see that the density function
of $\theta:=2a_2\lambda$ is
\begin{eqnarray*}
& & g_{\beta, a_1, a_2}(\theta)\\
 &:=& c_{J}^{\beta, a_1,
a_2}\left(\frac{1}{2a_2}\right)^{n(n-1)\beta/2+n(a_1-p)+n}\prod_{1\leq
i<j\leq
n}|\theta_i-\theta_j|^{\beta}\cdot\prod_{i=1}^n\theta_i^{a_1-p}
\left(1-\frac{\theta_i}{2a_2}\right)^{a_2-p}\\
& = & c_{J}^{\beta, a_1,
a_2}\left(\frac{1}{2a_2}\right)^{na_1}\prod_{1\leq
i<j\leq
n}|\theta_i-\theta_j|^{\beta}\cdot\prod_{i=1}^n\theta_i^{a_1-p}
\left(1-\frac{\theta_i}{2a_2}\right)^{a_2-p}
\end{eqnarray*}
for $0\leq \theta_i\leq 2a_2$ and $1\leq i \leq n$, and is equal to
zero, otherwise. Therefore,
\begin{eqnarray}
\|\ml{L}(2a_2\lambda)-\ml{L}(\mu)\| &= &
\int_{[0,\infty)^n}|g_{\beta,
a_1, a_2}(\mu)- f_{\beta, a_1}(\mu)|\,d\mu\nonumber\\
& = & \int_{[0,\infty)^n}|\frac{g_{\beta, a_1, a_2}(\mu)}{f_{\beta,
a_1}(\mu)}-1|\cdot f_{\beta, a_1}(\mu)\,d\mu.\lbl{lays}
\end{eqnarray}
Now, review $f_{\beta,
a_1}(\mu)$ as in (\ref{bWishart}) to see that
\begin{eqnarray*}
\frac{g_{\beta, a_1, a_2}(\mu)}{f_{\beta, a_1}(\mu)} =
\frac{c_{J}^{\beta, a_1, a_2}}{c_{L}^{\beta,
a_1}}\left(\frac{1}{2a_2}\right)^{na_1}\cdot\prod_{i=1}^n
\left(1-\frac{\mu_i}{2a_2}\right)^{a_2-p}\cdot
e^{\sum_{i=1}^n\mu_i/2}
\end{eqnarray*}
for $0\leq \mu_i\leq 2a_2,\ i=1,2,\cdots, n,$ and is zero,
otherwise. It is easy to check that
\begin{eqnarray*}
\frac{c_{J}^{\beta, a_1, a_2}}{c_{L}^{\beta,
a_1}}\left(\frac{1}{2a_2}\right)^{na_1}
& = & 2^{na_1}\cdot \left(\frac{1}{2a_2}\right)^{na_1}\prod_{j=1}^{n}\frac{\Gamma(a_1+a_2-\frac{\beta}{2}(n-j))}
{\Gamma(a_2-\frac{\beta}{2}(n-j))}\\
&= & a_2^{-na_1}\prod_{i=0}^{n-1}\frac{\Gamma(a_1+a_2-\frac{\beta}{2}i)}
{\Gamma(a_2-\frac{\beta}{2}i)}=K_n.
\end{eqnarray*}
Thus, $g_{\beta, a_1, a_2}(\mu)/f_{\beta, a_1}(\mu)=K_n\cdot L_n(\mu),$ which together with (\ref{celcius}) and (\ref{lays}) yields the first conclusion. Finally,
\begin{eqnarray*}
E(K_{n}\cdot L_n(\mu))=\int\frac{g_{\beta, a_1, a_2}(\mu)}{f_{\beta, a_1}(\mu)}\cdot f_{\beta, a_1}(\mu)\,d\mu =\int g_{\beta, a_1, a_2}(\mu)\,d\mu=1.\ \ \ \ \ \ \ \ \blacksquare
\end{eqnarray*}


\begin{lemma}\lbl{qian} Let $h(x)=x\log x$ for $x>0.$ For a fixed constant $\beta>0,$ an integer $n\geq 1$ and variables $a_1>0$ and $a_2>0$, set  $b_1=\frac{2}{\beta}a_1$ and $b_2=\frac{2}{\beta}a_2$. If $n\to \infty,\ a_1\to \infty$ and $a_2\to \infty $ in a way that
$a_1=o(\sqrt{a_2})$ and $n=o(\sqrt{a_2})$, then
\begin{eqnarray*}
\sum_{i=1}^n\left\{h(b_1+b_2-i+1)-h(b_2-i+1)\right\}=n b_1\left(1+\log  b_2+\frac{b_1-n}{2b_2}\right) +o(1).
\end{eqnarray*}
\end{lemma}
\textbf{Proof}. Note that  $h'(x)=1+\log x,\,
h''(x)=1/x$ and $h^{(3)}(x)=-1/x^2.$  Given $x_0>0,$ for any $\Delta x>-x_0,$  by the Taylor
expansion,
\begin{eqnarray*}
h(x_0+\Delta x)- h(x_0)& = & h'(x_0)\Delta x
+\frac{1}{2}h''(x_0)(\Delta
x)^2+ \frac{1}{6}h^{(3)}(\xi)(\Delta x)^3 \\
& = &  (1+\log x_0)\Delta x +\frac{1}{2x_0}(\Delta x)^2-
\frac{1}{6\xi^2}(\Delta x)^3
\end{eqnarray*}
where $\xi$ is between $x_0$ and $x_0+\Delta x$. Now take
$x_0=b_2-i+1$ and $\Delta x=b_1,$ we have that
\begin{eqnarray}
& & h(b_1+b_2-i+1)-h(b_2-i+1)\nonumber\\
&= & b_1(1+\log (b_2-i+1))+\frac{b_1^2}{2}\cdot\frac{1}{b_2-i+1} +
O\left(\frac{b_1^3}{b_2^2}\right)\lbl{lovegod}
\end{eqnarray}
uniformly for all $1\leq i \leq n$. Obviously,
\begin{eqnarray}\lbl{cooky}
\log (b_2-i+1)=\log b_2+
\log \left(1-\frac{i-1}{b_2}\right)=\log b_2 -\frac{i-1}{b_2} +
O\left(\frac{n^2}{b_2^2}\right)
\end{eqnarray}
uniformly over all $1\leq i \leq n$ as
\begin{eqnarray}\lbl{take}
n\to \infty,\ a_1\to \infty\ \mbox{and}\ a_2\to \infty\ \mbox{such that}\ a_1=o(\sqrt{a_2}),\mbox{and}\ n=o(\sqrt{a_2}).
\end{eqnarray}
Now,
\begin{eqnarray*}
\frac{b_1^2}{2}\cdot\frac{1}{b_2-i+1} & = & \frac{b_1^2}{2b_2} +
\frac{b_1^2}{2}\cdot\left(\frac{1}{b_2-i+1}-\frac{1}{b_2}\right)\\
& = & \frac{b_1^2}{2b_2} +
\frac{b_1^2}{2}\cdot\frac{i-1}{b_2(b_2-i+1)}\\
& = & \frac{b_1^2}{2b_2} + O\left(\frac{nb_1^2}{b_2^2}\right)
\end{eqnarray*}
uniformly for all $1\leq i \leq n$ as (\ref{take}) holds. Therefore, by (\ref{lovegod}) and (\ref{cooky}),
\begin{eqnarray*}
& & h(b_1+b_2-i+1)-h(b_2-i+1)\\
&= & b_1+ b_1\log b_2-\frac{b_1(i-1)}{b_2} + \frac{b_1^2}{2b_2}+
O\left(\frac{b_1^3 + n^2b_1 +  nb_1^2}{b_2^2}\right)
\end{eqnarray*}
uniformly for all $1\leq i \leq n$ as (\ref{take}) holds.  Thus
\begin{eqnarray*}
& & \sum_{i=1}^n\left\{h(b_1+b_2-i+1)-h(b_2-i+1)\right\}\\
&=  & nb_1 + nb_1\log b_2 - \frac{b_1n(n-1)}{2b_2} + \frac{nb_1^2}{2b_2} + n \cdot O\left(\frac{b_1^3 + n^2b_1 +  nb_1^2}{b_2^2}\right)\\
& = & nb_1\left(1+ \log b_2+\frac{b_1-n}{2b_2}\right) + \frac{b_1n}{2b_2} + n \cdot O\left(\frac{b_1^3 + n^2b_1 +  nb_1^2}{b_2^2}\right).
\end{eqnarray*}
The conclusion follows because the last two terms are all of order $o(1)$ as (\ref{take}) holds.\ \ \ \ \ \ \ \ $\blacksquare$

\begin{lemma}\lbl{ship} Let $K_n$ be as in (\ref{shining1}).  Assuming (\ref{shanxi}), we then have
\begin{eqnarray}\lbl{baby1}
K_n=\exp\left\{\frac{(1-\gamma)\beta^2n^3}{8a_2\gamma^2} + o(1)\right\}.
\end{eqnarray}

\end{lemma}
\textbf{Proof}. We claim that it suffices to prove
\begin{eqnarray}\lbl{baby2}
K_n=\exp\left\{\frac{na_1(a_1-\frac{\beta}{2} n)}{2a_2} + o(1)
\right\}
\end{eqnarray}
under assumption (\ref{take}).   If this is true, under the condition that $n\beta/(2a_1)\to \gamma\leq 1$, it is easy to check that
\begin{eqnarray*}
\frac{na_1(a_1-\frac{\beta}{2} n)}{2a_2} =\frac{\beta^2n^3t_n(t_n-1)}{8a_2}= \frac{(1-\gamma)\beta^2n^3}{8a_2\gamma^2} + o(1)
\end{eqnarray*}
as (\ref{take}) holds, where $t_n:=2a_1/(\beta n) \to \gamma^{-1}$. Thus, (\ref{baby1}) is obtained.

Now we prove (\ref{baby2}). Set $\alpha=\frac{\beta}{2},\
b_1=\frac{2}{\beta}a_1$ and $b_2=\frac{2}{\beta}a_2.$ It is easy to
see that
\begin{eqnarray*}
K_n=\left(\frac{1}{\alpha
b_2}\right)^{n\alpha b_1}\cdot\prod_{i=1}^{n}\frac{\Gamma(\alpha(b_1+b_2-i+1))}
{\Gamma(\alpha(b_2-i+1))}.
\end{eqnarray*}
Recall the Stirling formula:
\begin{eqnarray*}
\log\Gamma(z)=z\log z - z -\frac{1}{2}\log z+ \log \sqrt{2\pi}
+\frac{1}{12z} +O\left(\frac{1}{x^3}\right)
\end{eqnarray*}
as $x=\mbox{Re}\,(z)\to +\infty,$ where
$\Gamma(z)=\int_{0}^{\infty}t^{z-1} e^{-t}\, dt$ with
$\mbox{Re}\,(z)>0$, see, e.g., p.368 from \cite{Gamelin} or (37) on
p.204 from \cite{Ahlfors}. It follows that
\begin{eqnarray}
& & \log K_n\lbl{piano}\\
&=& -\alpha n b_1 \log (\alpha b_2)\nonumber\\
& & +\sum_{i=1}^n\left\{\alpha(b_1+b_2-i+1)\log \alpha(b_1+b_2-i+1)
-
\alpha(b_2-i+1)\log \alpha(b_2-i+1)\right\}\lbl{xiaohai1}\ \ \ \ \lbl{wei1}\\
 & & -\alpha n b_1 -\frac{1}{2}\sum_{i=1}^n\log \frac{b_1+b_2-i+1}{b_2-i+1}
 + O\left(\frac{1}{b_2-n}\right)\lbl{xiaohai22}
\end{eqnarray}
as (\ref{take}) holds.

Now, write $(\alpha x)\log (\alpha x)=(\alpha \log \alpha) x + \alpha (x
\log x)$ and set $h(x)=x\log x$ for $x>0.$ Calculating the difference between the two terms in
 the sum of (\ref{wei1}), we know that the whole sum in (\ref{wei1}) is
identical to
\begin{eqnarray}\lbl{Bowen1}
& & \alpha(\log \alpha) nb_1
+\alpha\sum_{i=1}^n(h(b_1+b_2-i+1)-h(b_2-i+1))\nonumber\\
& = & \alpha(\log \alpha) nb_1
+\alpha n b_1\left(1+\log  b_2+\frac{b_1-n}{2b_2}\right) + o(1)\nonumber\\
& = & \alpha n b_1 + (\alpha n b_1)\log (\alpha b_2) + \alpha n b_1\cdot \frac{b_1-n}{2b_2} + o(1)
\end{eqnarray}
by Lemma \ref{qian}. From the fact that $\log (1+x) \leq x$ for any
$x\geq 0,$ we have
\begin{eqnarray*}\lbl{lib}
0< \sum_{i=1}^n\log \frac{b_1+b_2-i+1}{b_2-i+1}=\sum_{i=1}^n\log
\left(1+\frac{b_1}{b_2 - i+1}\right)\leq \frac{nb_1}{b_2-n}
\end{eqnarray*}
for any $b_2>n.$ Thus, the sum of the three terms in (\ref{xiaohai22}) is equal to $-\alpha n b_1  +
O\left(\frac{nb_1}{b_2}\right)$ as (\ref{take}) holds. Combining this and (\ref{piano})-(\ref{Bowen1}), we get
\begin{eqnarray*}
\log K_n =\alpha n b_1\cdot \frac{b_1-n}{2b_2} + o(1)=\frac{na_1(a_1-\frac{\beta}{2} n)}{2a_2} + o(1)
\end{eqnarray*}
as (\ref{take}) holds. This gives (\ref{baby2}).\ \ \ \ \ \ \ \ $\blacksquare$ \\

\noindent In the proofs next, we use $o_P(1)$ to denote a random variable that goes zero in probability as taking a limit.
\begin{lemma}\lbl{jiwang}  Let  $\mu=(\mu_1,\cdots,
\mu_n)$ be a random variable with the density function as in
(\ref{bWishart}) with $a=a_1$. Let $L_n(\mu)$ be as in (\ref{shining2}).
If (\ref{shanxi}) holds, then
\begin{eqnarray*}
\exp\left\{\frac{(1-\gamma)\beta^2n^3}{8a_2\gamma^2} \right\}\cdot L_n(\mu) \to 1
\end{eqnarray*}
in probability as $n\to \infty$.
\end{lemma}
\noindent\textbf{Proof}. From (\ref{shining2}), we see that
\begin{eqnarray*}
L_n(\mu)=e^{(1/2)\sum_{i=1}^n{\mu_i}}\cdot
\prod_{i=1}^n\left(1-\frac{\mu_i}{2a_2}\right)^{a_2-p}
 \cdot I(\max_{1\leq
i \leq n}\mu_i \leq 2a_2).
\end{eqnarray*}
By Lemma \ref{repair}, since $\frac{n\beta}{2a_1}\to \gamma\in (0, 1]$ by (\ref{shanxi}),
\begin{eqnarray}\lbl{berger}
\frac{\max_{1\leq
i \leq n}\mu_i}{n} \to \beta (1+\sqrt{\gamma^{-1}})^2
\end{eqnarray}
in probability as $n\to\infty.$ Since $n=o(\sqrt{a_2}\,),$ choose $\delta_n = (n\sqrt{a_2}\,)^{1/2}$, then $\delta_n/n \to \infty$ and $\delta_n/\sqrt{a_2}\to 0$ as taking the limit as in (\ref{shanxi}).   Therefore, to prove the lemma, it is enough to show
\begin{eqnarray}\lbl{wisconsin}
\exp\left\{\frac{(1-\gamma)\beta^2n^3}{8a_2\gamma^2} \right\}\cdot \tilde{L}_n(\mu) \to 1
\end{eqnarray}
in probability as $n\to\infty$, where
\begin{eqnarray}\lbl{light}
\tilde{L}_n(\mu):=e^{(1/2)\sum_{i=1}^n{\mu_i}}\cdot
\prod_{i=1}^n\left(1-\frac{\mu_i}{2a_2}\right)^{a_2-p}
 \cdot I(\max_{1\leq
i \leq n}\mu_i \leq \delta_n).
\end{eqnarray}
This is because, for any two sequences random variables $\{\xi_n;\, n\geq 1\}$ and $\{\eta_n;\, n\geq 1\}$, if $\xi_n \to 1$ in probability and $P(\xi_n \neq \eta_n)\to 0$ as $n\to\infty$, then $\eta_n\to 1$ in probability as $n\to\infty.$ Rewrite
\begin{eqnarray*}
 \tilde{L}_n(\mu) &
= & \exp\left\{\frac{1}{2}\sum_{i=1}^n{\mu_i} + (a_2-p)\sum_{i=1}^n\log (1-\frac{\mu_i}{2a_2})\right\}
\end{eqnarray*}
on $\Omega_n:=\{\max_{1\leq
i \leq n}\mu_i \leq \delta_n\}$.
 Noticing $\log (1-x)=-x-(x^2/2) +O(x^3)$ as $x\to 0$,
\begin{eqnarray}\lbl{temple}
\sum_{i=1}^n\log (1-\frac{\mu_i}{2a_2}) =-\frac{1}{2a_2}\sum_{i=1}^n\mu_i - \frac{1}{8a_2^2}\sum_{i=1}^n\mu_i^2 + O\left(\frac{1}{a_2^3}\sum_{i=1}^n\mu_i^3\right)
\end{eqnarray}
on $\Omega_n.$  Now, on $\Omega_n$ again,
\begin{eqnarray}\lbl{river5}
\frac{1}{a_2^3}\sum_{i=1}^n\mu_i^3 \leq  \frac{n(\delta_n)^3}{a_2^3}=\left(\frac{\delta_n}{\sqrt{a_2}}\right)^3\cdot \frac{n}{\sqrt{a_2}}\cdot \frac{1}{a_2}\to 0
\end{eqnarray}
as taking the limit in (\ref{shanxi}). Recall $p=1+\frac{\beta}{2}(n-1)$. We have from (\ref{temple}) and (\ref{river5}) that, on $\Omega_n$,
\begin{eqnarray*}
& & (a_2-p)\sum_{i=1}^n\log (1-\frac{\mu_i}{2a_2})\\
& = & -\frac{a_2-p}{2a_2}\left(-\frac{\beta n^2}{\gamma} + \sum_{i=1}^n\mu_i\right) - \frac{(a_2-p)\beta n^2}{2a_2\gamma} \\
& & -\frac{a_2-p}{8a_2^2}\left(-\frac{\beta^2n^3}{\gamma^2}(1+\gamma)+\sum_{i=1}^n\mu_i^2\right) - \frac{\beta^2n^3(a_2-p)}{8a_2^2\gamma^2}(1+\gamma) + O\left((\frac{\delta_n}{\sqrt{a_2}})^3\cdot \frac{n}{\sqrt{a_2}}\right)
\end{eqnarray*}
as (\ref{shanxi}) holds. By Lemma \ref{flu}, as $n\to \infty$,
\begin{eqnarray}
& & \frac{1}{n^2}\sum_{i=1}^n\mu_i \overset{P}{\to}
 \frac{\beta}{\gamma};\ \ \frac{1}{n}\sum_{i=1}^n\mu_i- \frac{\beta}{\gamma}n \Rightarrow N(0, \sigma_1^2)\lbl{useful1};\\
& & \frac{1}{n^3}\sum_{i=1}^n\mu_i^2 \overset{P}{\to} \frac{\beta^2}{\gamma^2}(1+\frac{1}{2}\cdot 2\gamma)=\frac{\beta^2}{\gamma^2}(1+\gamma);\ \ \frac{1}{n^2}\sum_{i=1}^n\mu_i^2- \frac{\beta^2}{\gamma^2}(1+\gamma)n\Rightarrow N(0, \sigma_2^2)\ \ \ \ \ \ \ \ \lbl{useful2}
\end{eqnarray}
where $\sigma_1, \sigma_2$ are constants depending on $\gamma$ only,  the notation ``$\overset{P}{\to}$" means  ``converges in probability to" and ``$\Rightarrow$" means ``converges weakly to". Now, write $(a_2-p)/2a_2=(1/2) -p/2a_2, $ then
\begin{eqnarray*}
& & -\frac{a_2-p}{2a_2}\left(-\frac{\beta n^2}{\gamma} + \sum_{i=1}^n\mu_i\right)\\
& = & \frac{\beta n^2}{2\gamma}-\frac{1}{2}\sum_{i=1}^n\mu_i+\frac{pn}{2a_2}\cdot \frac{1}{n}\left(-\frac{\beta n^2}{\gamma} + \sum_{i=1}^n\mu_i\right)\\
& = & \frac{\beta n^2}{2\gamma}-\frac{1}{2}\sum_{i=1}^n\mu_i+ o_P(1)
\end{eqnarray*}
by (\ref{useful1}) as (\ref{shanxi}) holds. Also, under the same condition,  $(a_2-p)n^2/a_2^2=O(n^2/a_2)=o(1)$. It follows from  (\ref{useful2}) that
\begin{eqnarray*}
& & -\frac{a_2-p}{8a_2^2}\left(-\frac{\beta^2n^3}{\gamma^2}(1+\gamma)+\sum_{i=1}^n\mu_i^2\right)\\
&=& -\frac{(a_2-p)n^2}{8a_2^2}\left(-\frac{\beta^2n}{\gamma^2}(1+\gamma)+\frac{1}{n^2}\sum_{i=1}^n\mu_i^2\right)\to 0
\end{eqnarray*}
in probability as taking the limit in (\ref{shanxi}). In summary, combining all the computations above,
\begin{eqnarray*}
 & & \frac{1}{2}\sum_{i=1}^n\mu_i + (a_2-p)\sum_{i=1}^n\log (1-\frac{\mu_i}{2a_2})\\
 & = & \frac{\beta n^2}{2\gamma}-\frac{(a_2-p)\beta n^2}{2a_2\gamma}-\frac{\beta^2n^3(a_2-p)}{8a_2^2\gamma^2}(1+\gamma) + o_P(1)\\
 & = &  \frac{\beta pn^2}{2a_2\gamma} - \frac{\beta^2n^3}{8a_2\gamma^2}(1+\gamma) + \frac{\beta^2n^3p}{8a_2^2\gamma^2}(1+\gamma)  + o_P(1)
\end{eqnarray*}
on $\Omega_n$. Now, since $p=1+\frac{\beta}{2}(n-1)$, $n/\sqrt{a_2} \to 0$,  we see that
\begin{eqnarray*}
\frac{\beta pn^2}{2a_2\gamma}=\frac{\beta^2n^3}{4a_2\gamma} + o(1)\ \ \ \mbox{and}\ \ \ \frac{\beta^2n^3p}{8a_2^2\gamma^2}(1+\gamma) \to 0
\end{eqnarray*}
as (\ref{shanxi}) holds. Thus, on $\Omega_n$,
\begin{eqnarray*}
& & \frac{1}{2}\sum_{i=1}^n\mu_i +(a_2-p)\sum_{i=1}^n\log (1-\frac{\mu_i}{2a_2})\\
&= & \frac{\beta^2n^3}{4a_2\gamma}-\frac{\beta^2n^3}{8a_2\gamma^2}(1+\gamma) + o_P(1)=\frac{(\gamma-1)\beta^2n^3}{8a_2\gamma^2} + o_P(1)
\end{eqnarray*}
 as taking the limit in (\ref{shanxi}). By reviewing (\ref{light}), we conclude (\ref{wisconsin}).\ \ \ \ \ \ \ $\blacksquare$\\

The following is a variant of the Scheffe Lemma .
\begin{lemma}\lbl{cutie} Let $\{X_n;\, n\geq 1\}$ be a sequence of non-negative random variables. If $X_n \to 1$  in probability and $EX_n\to 1$ as $n\to \infty$, then $E|X_n-1| \to 0$ as $n\to\infty.$
\end{lemma}
\textbf{Proof}.  By the Skorohod representation theorem (see, e.g., p.85 from \cite{Durrett}), w.l.o.g., assume $X_n\to 1$ almost surely as $n\to\infty$. Thus, for any $K>1,$ we have $X_nI(X_n\leq K) \to 1$ almost surely as $n\to\infty$. This gives that $E|X_nI(X_n\leq K) - 1|\to 0$ and $EX_nI(X_n\leq K) \to 1$ as $n\to\infty.$ The second assertion and the condition that $EX_n \to 1$ imply  $EX_nI(X_n> K) \to 0$ as $n\to\infty.$ Therefore,
\begin{eqnarray*}
E|X_n-1|\leq E|X_nI(X_n\leq K) - 1| + EX_nI(X_n> K) \to 0
\end{eqnarray*}
as $n\to\infty$.\ \ \ \ \ \ \ \ \ \ $\blacksquare$\\

\noindent\textbf{Proof of Theorem \ref{main}}. It is known from Lemma \ref{botany} that
\begin{eqnarray*}
\|\ml{L}(2a_2\lambda)-\ml{L}(\mu)\|=E|K_{n}\cdot
L_n(\mu)-1|.
\end{eqnarray*}
with $E(K_n\cdot L_n(\mu))=1$ for all $n\geq 2,$ where $\mu$ has density $f_{\beta,a_1}(\mu)$ as in (\ref{bWishart}). By Lemmas \ref{ship} and \ref{jiwang},
\begin{eqnarray*}
K_n=\exp\left\{\frac{(1-\gamma)\beta^2n^3}{8a_2\gamma^2} + o(1)\right\}\ \ \mbox{and} \ \ \ \exp\left\{\frac{(1-\gamma)\beta^2n^3}{8a_2\gamma^2} \right\}\cdot L_n(\mu) \to 1
\end{eqnarray*}
in probability as taking the limit in (\ref{shanxi}). These imply that $K_{n}\cdot
L_n(\mu)\to 1$ in probability as taking the same limit. Then the desired conclusion follows from Lemma \ref{cutie}.\ \ \ \ \ \ $\blacksquare$\\

\section{The Proofs of Main Results}\lbl{shortproof}
\setcounter{equation}{0}
By using Theorem \ref{main} developed in Section \ref{yang}, we now are ready to prove the results stated in  Introduction.

\noindent\textbf{Proof of Theorem \ref{snowf}}. Set
\begin{eqnarray*}
\nu_n=\frac{1}{n}\sum_{i=1}^nI_{\frac{a_2}{n}\lambda_i'}
\end{eqnarray*}
for $(\lambda_1', \cdots, \lambda_n')\in [0,+\infty)^n.$ Then, recall the definition of $d$ in (\ref{du}), by the triangle inequality,
\begin{eqnarray*}
|d(\mu_n, \mu_0) - d(\nu_n, \mu_0)| & \leq & d(\mu_n, \nu_n)\\
&=& \sup_{\|f\|_{BL}\leq 1}|\frac{1}{n}\sum_{i=1}^n\left(f(n^{-1}a_2\lambda_i)-f(n^{-1}a_2\lambda_i')\right)|\\
& \leq & \frac{a_2}{n}\cdot \max_{1\leq i \leq n}|\lambda_i-\lambda_i'|,
\end{eqnarray*}
where the Lipschitz inequality $|f(x)-f(y)|\leq |x-y|$ is used in the last step. This says that $d(\mu_n, \mu_0)$, as a function of $(\lambda_1, \cdots \lambda_n)$, is continuous for each $n\geq 2$. Thus, for any $\epsilon>0,$ there exists a (non-random) Borel set $A\subset \mathbb{R}^n$ such that $\{d(\mu_n, \mu_0)\geq \epsilon\}=\{(\lambda_1, \cdots, \lambda_n)\in A\}$. Then, by the definition of the variation norm in  (\ref{variation}) we see that
\begin{eqnarray}\lbl{brown}
P(d(\mu_n, \mu_0)\geq \epsilon) \leq P(d(\mu_n', \mu_0)\geq \epsilon) + \|\mathcal{L}(2a_2\lambda) - \mathcal{L}(\mu)\|
\end{eqnarray}
for any $\epsilon>0,$ where $\mu_n'=(1/n)\sum_{i=1}^nI_{\mu_i/(2n)}$ and $\mu=(\mu_1,\cdots, \mu_n)$ has density $f_{\beta, a_1}(\mu)$ as in (\ref{bWishart}) with $a=a_1$ and $n\beta/2a_1 \to \gamma \in (0, 1]$. By Lemma \ref{xiaochild}, with probability one,
\begin{eqnarray}\lbl{violin}
\frac{1}{n}\sum_{i=1}^nI_{\frac{\mu_i\gamma}{n\beta}}\ \ \mbox{converges weakly to}\ \ \mu_{\infty}
\end{eqnarray}
with density $f_{\gamma}(x)$ as in (\ref{mooncake}). Write $\mu_i/(2n)=(\mu_i\gamma/n\beta)c^{-1}$, where $c=2\gamma/\beta$. Then, by (\ref{violin}), with probability one, $\mu_n'$ converges weakly to $\mu_0$, where $\mu_0$ has density function $c\cdot f_{\gamma}(cx).$ Equivalently, $d(\mu_n', \mu_0) \to 0$ almost surely. This, (\ref{brown}) and Theorem \ref{main} in Section \ref{yang} prove the theorem.\ \ \ \ \ \ \ $\blacksquare$\\

\noindent\textbf{Proof of Theorem \ref{branch}}. First, $\lambda_{max}(n)$ and $\lambda_{min}(n)$ are continuous functions of $\lambda_1, \cdots, \lambda_n$ for any $n\geq 1.$ Then
\begin{eqnarray}\lbl{greensnow}
& & P\left(|\frac{a_2}{n}\lambda_{max}(n) - \frac{\beta (1+\sqrt{\gamma})^2}{2\gamma}|\geq \epsilon\right)\nonumber\\
& \leq & P\left(|\frac{1}{2n}\mu_{max}(n) - \frac{\beta (1+\sqrt{\gamma})^2}{2\gamma}|\geq \epsilon\right)  + \|\mathcal{L}(2a_2\lambda) - \mathcal{L}(\mu)\|
\end{eqnarray}
for any $\epsilon>0,$ where  $\mu=(\mu_1,\cdots, \mu_n)$ has density $f_{\beta, a_1}(\mu)$ as in (\ref{bWishart}) with $a=a_1$ and $n\beta/2a_1 \to \gamma \in (0, 1]$. From Lemma \ref{repair}, we know $\mu_{\max}(n)/(2n)\to \beta (1+\sqrt{\gamma^{-1}})^2/2=\beta (1+\sqrt{\gamma})^2/(2\gamma)$ in probability. This together with (\ref{greensnow}) and Theorem \ref{main} in Section \ref{yang} yields the desired conclusion. By the same argument, $(a_2/n)\lambda_{min}(n)$ converges to $\beta(1-\sqrt{\gamma})^2/(2\gamma)$ in probability.\ \ \ \ \ \ \ $\blacksquare$\\

\noindent\textbf{Proof of Theorem \ref{hay}}. Evidently,
\begin{eqnarray*}
\left|P((X_1,\cdots, X_k)\in A)- P((Y_1,\cdots, Y_k)\in A)\right| \leq \|\mathcal{L}(2a_2\lambda) - \mathcal{L}(\mu)\|
\end{eqnarray*}
for any Borel set $A \in \mathbb{R}^k$, where
\begin{eqnarray*}
Y_i
&= & \sum_{j=1}^n\left(\frac{c}{2n}\mu_j\right)^i-n
\sum_{r=0}^{i-1}\frac{1}{r+1}\binom{i}{r}\binom{i-1}{r}\gamma^{r}\\
&= & \sum_{j=1}^n\left(\frac{\gamma}{n\beta}\mu_j\right)^i-n
\sum_{r=0}^{i-1}\frac{1}{r+1}\binom{i}{r}\binom{i-1}{r}\gamma^{r}
\end{eqnarray*}
for $i\geq 1$ (since $c=2\gamma/\beta$), and  $\mu=(\mu_1,\cdots, \mu_n)$ has density $f_{\beta, a_1}(\mu)$ as in (\ref{bWishart}) with $a=a_1$ and $n\beta/2a_1 \to \gamma \in (0, 1]$. The conclusion then follows from this, Theorem A.\ref{haci} and Theorem \ref{main} in Section \ref{yang}.\ \ \ \ \ \ \ $\blacksquare$\\

\noindent\textbf{Proof of Theorem \ref{longlive}.} The assumption that $2\beta^{-1}a_1-n=c$ and $n=o(\sqrt{a_2})$ imply  that  $n\beta/2a_1 \to 1$ and $a_1=o(\sqrt{a_2}).$ Thus Theorem \ref{main} in Section \ref{yang} holds.

Let $(\theta_1, \cdots, \theta_n)$ have density $f_{\beta, a_1}$ as in (\ref{bWishart}) with $a=a_1$. Noticing,  ``$\beta \lambda_i$" and ``$a$" in Theorem 1 from \cite{Rider} correspond to ``$\theta_i$" and ``$c$" here, respectively. By Theorem 1 from \cite{Rider}, for fixed integer  $k\geq 1$,
\begin{eqnarray*}
\left(\frac{n}{\beta}\theta^{(n)}, \cdots, \frac{n}{\beta}\theta^{(n-k+1)}\right)\ \ \mbox{converges weakly to}\ \ \ (\Lambda_0(\beta, c), \cdots,\Lambda_{k-1}(\beta, c) )
\end{eqnarray*}
as $n\to\infty.$ By Theorem \ref{main} in Section \ref{yang},
\begin{eqnarray*}
P\left((2a_2\lambda_1, \cdots, 2a_2\lambda_n)\in B_n\right) - P((\theta_1, \cdots, \theta_n)\in B_n)\to 0
\end{eqnarray*}
for any Borel set $B_n\subset \mathbb{R}^n$ for $n\geq 1.$ From the Weyl perturbation theorem, we know that $\lambda^{(i)}$ is a continuous function of $(\lambda_1, \cdots, \lambda_n)$ for any $1\leq i \leq n$. Combining the above two limits we obtain
\begin{eqnarray*}
\left(\frac{2a_2n}{\beta}\lambda^{(n)}, \cdots, \frac{2a_2n}{\beta}\lambda^{(n-k+1)}\right)\ \ \mbox{converges weakly to}\ \ \ (\Lambda_0(\beta, c), \cdots,\Lambda_{k-1}(\beta, c) )
\end{eqnarray*}
as $n\to\infty$ and $a_2 \to\infty$ with $n=o(\sqrt{a_2}).$ The proof is complete.\ \ \ \ \ \ \ \ $\blacksquare$\\

\noindent\textbf{Proof of Theorem \ref{work}}. Recalling (\ref{bWishart}), let
\begin{eqnarray}\lbl{attack}
\tilde{m}_n=\left(\sqrt{n} +\sqrt{2\beta^{-1}a}\,\right)^2\ \ \mbox{and}\ \ \tilde{\sigma}_n=\frac{(2\beta^{-1}na)^{1/6}}{(\sqrt{n} +\sqrt{2\beta^{-1}a}\ )^{4/3}}.
\end{eqnarray}
Let $(\theta_1, \cdots, \theta_n)$ have density $f_{\beta, a}$ as in (\ref{bWishart}). Noticing,  ``$\beta \lambda_i$" in Theorem 1.4 from \cite{RRV} corresponds to ``$\theta_i$" here; ``$\kappa$" in Theorem 1.4 from \cite{RRV} is equal to $2\beta^{-1}a$, and $\frac{\beta}{2}(n-1) + 1=p,$ and the $k$-th lowest eigenvalue of $\mathcal{H}_{\beta}$ is the $(n-k+1)$-th largest eigenvalue of $-\mathcal{H}_{\beta}$.  Then by Theorem 1.4 from \cite{RRV},
\begin{eqnarray}\lbl{clark}
\tilde{\sigma}_n \left(\frac{\theta^{(l)}}{\beta}-\tilde{m}_n\right)_{l=1,\cdots,k}\ \ \mbox{converges weakly to}\ \ (\Lambda_1, \cdots, \Lambda_k)
\end{eqnarray}
as $n\to\infty$ and $a\to\infty$ such that $n/a$ converges to a nonzero, finite constant. In other words,
\begin{eqnarray*}
P\left(\tilde{\sigma}_n \left(\frac{\theta^{(l)}}{\beta}-\tilde{m}_n\right)_{l=1,\cdots,k} \in A\right) \to P\left((\Lambda_1, \cdots, \Lambda_k)\in A\right)
\end{eqnarray*}
for any Borel set $A \subset \mathbb{R}^k.$  By Theorem \ref{main} in Section \ref{yang}, assuming (\ref{shanxi}) and $a=a_1$,
\begin{eqnarray*}
P\left((2a_2\lambda_1, \cdots, 2a_2\lambda_n)\in B_n\right) - P((\theta_1, \cdots, \theta_n)\in B_n)\to 0
\end{eqnarray*}
for any Borel set $B_n\subset \mathbb{R}^n$ for $n\geq 1.$ The Weyl perturbation theorem says that $g(x):=x^{(l)}$,  the $l$-th largest one in $\{x_1, \cdots, x_m\}$, is a continuous function of $(x_1, \cdots, x_m)\in \mathbb{R}^m$ for any $1\leq l\leq m$. Replacing $a$ by $a_1$ in (\ref{attack}), the above two assertions conclude that
\begin{eqnarray*}
\sigma_n \left(\frac{2a_2\lambda^{(l)}}{\beta}-m_n\right)_{l=1,\cdots,k}\ \ \mbox{converges weakly to}\ \ (\Lambda_1, \cdots, \Lambda_k).\ \ \ \ \ \ \ \ \ \ \ \ \blacksquare\\
\end{eqnarray*}

\section{Appendix}\lbl{appendix}
\setcounter{equation}{0}
Let $Q(x)$ be a real continuous function on $[0, \infty)$ such that for any $\epsilon>0$
\begin{eqnarray}\lbl{LDP}
\lim_{x\to +\infty} x e^{-\epsilon Q(x)}=0.
\end{eqnarray}
For each integer $n\geq 1,$ let $p(n)$ be an integer depending on $n.$ Let $\lambda_1, \cdots, \lambda_{p(n)}$ be non-negative random variables with joint probability density
\begin{eqnarray}
\nu_n:=\frac{1}{Z_n}\exp\left(-n\sum_{i=1}^{p(n)}Q(\lambda_i)\right)\prod_{i=1}^{p(n)}\lambda_i^{\gamma(n)}\prod_{1\leq i < j\leq p(n)}|\lambda_i-\lambda_j|^{\beta}
\end{eqnarray}
where $\beta>0$ is fixed and $\gamma(n)\geq 0$ depends on $n$. Let $\mu_n$ be the empirical probability measure of $\lambda_1, \cdots, \lambda_{p(n)}.$ Under the weak topology, the large deviations for $\{\mu_n\}$ is given below. For a reference of general large deviations, one can see, e.g., \cite{DZ98} and \cite{DW}.

\begin{theoremm}(Theorem 1 in \cite{Hiai})\lbl{casino} Assume $p(n)/n\to \kappa\in (0, 1]$ and $\gamma(n)/n\to \tau\geq 0$ as $n\to\infty.$ Then the finite limit $B:=\lim_{n\to\infty}n^{-2}\log Z_n$ exists and $\{\mu_n;\ n\geq 1\}$ satisfies the large deviation principle with speed $\{n^{-2};\ n\geq 1\}$ and good rate function
\begin{eqnarray}\lbl{rate}
I(\nu):=-\frac{\kappa^2\beta}{2}\iint\log |x-y|\, d\nu(x)\,d\nu(y) + \kappa\int(Q(x)-\tau\log x)\,d\nu(x) +B
\end{eqnarray}
for all probability measure $\nu$ defined on $[0, \infty)$. Moreover, there exists a unique probability measure $\nu_0$ on $[0, \infty)$ such that $I(\nu_0)=0.$
\end{theoremm}
In Theorem 1 from \cite{Hiai} or Theorem 5.5.1 from \cite{Hiaipetz}, the limit $\tau$ above is required to be strictly positive. However, after a check, it is found that the conclusion also holds for $\tau=0.$
\begin{theoremm} (part of Theorem 8 from \cite{Hiai})\lbl{ok} Let $\gamma \in (0,1]$,  $\gamma_{min}=(1-\sqrt{\gamma})^2$ and $\gamma_{max}=(1+\sqrt{\gamma})^2$. For probability measure $\nu$, define
\begin{eqnarray*}
J(\nu)=-\frac{\gamma^2}{2}\iint\log |x-y|\, d\nu(x)d\nu(y) + \frac{\gamma}{2}\int\left(x-(1-\gamma)\log x\right)d\nu(x).
\end{eqnarray*}
Then the unique minimizer of $J(\nu)$ over all probability measures on $[0, +\infty)$ is the Marchenko-Pastur law with  density function $f_{\gamma}(x)$ as in (\ref{mooncake}).
\end{theoremm}
\textbf{Proof}. Take $a$ in \cite{Hiai} equal to $\gamma$. Notice
\begin{eqnarray*}
\sqrt{4\gamma -(x-1-\gamma)^2}=\sqrt{(x-(1-\sqrt{\gamma})^2))((1+\sqrt{\gamma})^2-x)}=\sqrt{(x-\gamma_{min})(\gamma_{max}-x)}
\end{eqnarray*}
for all $x\in [\gamma_{min}, \gamma_{max}].$ Also, over all probability measure $\nu$ on $[0, \infty),$ taking the minimum for
\begin{eqnarray*}
I(\nu): &= & -\frac{\gamma^2}{2}\int\,\int\log |x-y|\, d\nu(x)\,d\nu(y) + \frac{\gamma}{2}\int(x-(1-\gamma)\log x)\,d\nu(x)\\
& & \ \ \ \ \ \ \ \ \ \ \ -\frac{1}{4}\left(3\gamma-\gamma^2\log \gamma +(1-\gamma)^2\log (1-\gamma)\right)
\end{eqnarray*}
is the same as doing so for $J(\nu),$ where $0\log 0:=0$ as the convention.  Then the conclusion follows from Theorem 8 in  \cite{Hiai}.\ \ \ \ \ \ \ \ $\blacksquare$

\begin{theoremm}\lbl{haci}(Theorem 1.5 from \cite{DumitriuGlobal}) Let  $\lambda_1,\cdots,
\lambda_n$ be random variables with the density function as in
(\ref{bWishart}). Assume $n\beta/(2a)\to \gamma\leq 1.$
Define
\begin{eqnarray*}
X_i=\sum_{j=1}^n\left(\frac{\gamma}{n\beta}\lambda_j\right)^i-n
\sum_{r=0}^{i-1}\frac{1}{r+1}\binom{i}{r}\binom{i-1}{r}\gamma^{r}
\end{eqnarray*}
for $i\geq 1.$ Then, as $n\to\infty$, $(X_1, \cdots, X_k)$ converges weakly to a normal distribution $N_k(\mu, \Sigma)$ for some $\mu$ and $\Sigma$.
\end{theoremm}

\noindent\textbf{Acknowledgement}\ The author thanks Professor Ioana Dumitriu for helpful discussions on the beta ensembles.\\


\baselineskip 10pt
\def\ref{\par\noindent\hangindent 30pt}

\end{document}